\documentclass[final,a4paper,11pt]{article}

\usepackage{xspace}
\usepackage{graphics} 
\usepackage{epsfig}
\usepackage{appendix}

\usepackage{floatflt}

\newcommand{\tlowername}[2]%
{$\stackrel{\makebox[1pt]{#1}}%
{\begin{picture}(0,0)%
\put(0,0){\makebox(0,6)[t]{\makebox[1pt]{$#2$}}}%
\end{picture}}$}%

\newcommand{\Tcase}[2]{\makebox[23pt]%
{\raisebox{2.5pt}{$\stackrel{#2}{#1{20}}$}}}%

\newcommand{\AR}[1]%
{\begin{picture}(#1,0)%
\put(0,0){\vector(1,0){#1}}%
\end{picture}}%

\newcommand{\DOTAR}[1]%
{\NUMBEROFDOTS=#1%
\divide\NUMBEROFDOTS by 3%
\begin{picture}(#1,0)%
\multiput(0,0)(3,0){\NUMBEROFDOTS}{\circle*{1}}%
\put(#1,0){\vector(1,0){0}}%
\end{picture}}%

\newcommand{\MONO}[1]%
{\begin{picture}(#1,0)%
\put(0,0){\vector(1,0){#1}}%
\put(2,-2){\line(0,1){4}}%
\end{picture}}%

\newcommand{\EPI}[1]%
{\begin{picture}(#1,0)(-#1,0)%
\put(-#1,0){\vector(1,0){#1}}%
\put(-6,-2){\line(0,1){4}}%
\end{picture}}%

\newcommand{\BIMO}[1]%
{\begin{picture}(#1,0)(-#1,0)%
\put(-#1,0){\vector(1,0){#1}}%
\put(-6,-2){\line(0,1){4}}%
\put(-#1,-2){\hspace{2pt}\line(0,1){4}}%
\end{picture}}%

\newcommand{\BIAR}[1]%
{\begin{picture}(#1,4)%
\put(0,0){\vector(1,0){#1}}%
\put(0,4){\vector(1,0){#1}}%
\end{picture}}%

\newcommand{\EQL}[1]%
{\begin{picture}(#1,0)%
\put(0,1){\line(1,0){#1}}%
\put(0,-1){\line(1,0){#1}}%
\end{picture}}%

\newcommand{\ADJAR}[1]%
{\begin{picture}(#1,4)%
\put(0,0){\vector(1,0){#1}}%
\put(#1,4){\vector(-1,0){#1}}%
\end{picture}}%

\newcommand{\TRIAR}[4]%
{%
\begin{tabular}{c}%
\mbox{$\scriptstyle#2$}\\ \AR{#1}\\%
\mbox{$\scriptstyle#3$}\\ \AR{#1}\\%
\mbox{$\scriptstyle#4$}\\ \AR{#1}%
\end{tabular}%
}%

\newcommand{\TRIADJAR}[4]%
{%
\begin{tabular}{c}%
\mbox{$\scriptstyle#2$}\\ \BKAR{#1}\\%
\mbox{$\scriptstyle#3$}\\ \AR{#1}\\%
\mbox{$\scriptstyle#4$}\\ \BKAR{#1}%
\end{tabular}%
}%

\newcommand{\QUADRIAR}[1]%
{%
\begin{tabular}{c}%
\AR{#1}\\%
\rule{0pt}{3pt}\\ \AR{#1}\\%
\rule{0pt}{3pt}\\ \AR{#1}\\%
\rule{0pt}{3pt}\\ \AR{#1}%
\end{tabular}%
}%

\newcommand{\QUADRIADJAR}[1]%
{%
\begin{tabular}{c}%
\BKAR{#1}\\%
\rule{0pt}{3pt}\\ \AR{#1}\\%
\rule{0pt}{3pt}\\ \BKAR{#1}\\%
\rule{0pt}{3pt}\\ \AR{#1}%
\end{tabular}%
}%

\newcommand{\QUINTIAR}[1]%
{%
\begin{tabular}{c}%
\AR{#1}\\%
\rule{0pt}{2pt}\\ \AR{#1}\\%
\rule{0pt}{2pt}\\ \AR{#1}\\%
\rule{0pt}{2pt}\\ \AR{#1}\\%
\rule{0pt}{2pt}\\ \AR{#1}%
\end{tabular}%
}%

\newcommand{\QUINTIADJAR}[1]%
{%
\begin{tabular}{c}%
\AR{#1}\\%
\rule{0pt}{2pt}\\ \BKAR{#1}\\%
\rule{0pt}{2pt}\\ \AR{#1}\\%
\rule{0pt}{2pt}\\ \BKAR{#1}\\%
\rule{0pt}{2pt}\\ \AR{#1}%
\end{tabular}%
}%

\newcommand{\Ar}[1]{\Tcase{\AR}{#1}}%

\newcommand{\BKAR}[1]%
{\begin{picture}(#1,0)%
\put(#1,0){\vector(-1,0){#1}}%
\end{picture}}%

\newcommand{\BKDOTAR}[1]%
{\NUMBEROFDOTS=#1%
\divide\NUMBEROFDOTS by 3%
\begin{picture}(#1,0)%
\multiput(#1,0)(-3,0){\NUMBEROFDOTS}{\circle*{1}}%
\put(0,0){\vector(-1,0){0}}%
\end{picture}}%

\newcommand{\BKMONO}[1]%
{\begin{picture}(#1,0)(-#1,0)%
\put(0,0){\vector(-1,0){#1}}%
\put(-2,-2){\line(0,1){4}}%
\end{picture}}%

\newcommand{\BKEPI}[1]%
{\begin{picture}(#1,0)%
\put(#1,0){\vector(-1,0){#1}}%
\put(6,-2){\line(0,1){4}}%
\end{picture}}%

\newcommand{\BKBIMO}[1]%
{\begin{picture}(#1,0)%
\put(#1,0){\vector(-1,0){#1}}%
\put(6,-2){\line(0,1){4}}%
\put(#1,-2){\hspace{-2pt}\line(0,1){4}}%
\end{picture}}%

\newcommand{\BKBIAR}[1]%
{\begin{picture}(#1,4)%
\put(#1,0){\vector(-1,0){#1}}%
\put(#1,4){\vector(-1,0){#1}}%
\end{picture}}%

\newcommand{\BKADJAR}[1]%
{\begin{picture}(#1,4)%
\put(0,4){\vector(1,0){#1}}%
\put(#1,0){\vector(-1,0){#1}}%
\end{picture}}%

\newcommand{\BKTRIAR}[4]%
{%
\begin{tabular}{c}%
\mbox{$\scriptstyle#2$}\\ \BKAR{#1}\\%
\mbox{$\scriptstyle#3$}\\ \BKAR{#1}\\%
\mbox{$\scriptstyle#4$}\\ \BKAR{#1}%
\end{tabular}%
}%

\newcommand{\BKTRIADJAR}[4]%
{%
\begin{tabular}{c}%
\mbox{$\scriptstyle#2$}\\ \AR{#1}\\%
\mbox{$\scriptstyle#3$}\\ \BKAR{#1}\\%
\mbox{$\scriptstyle#4$}\\ \AR{#1}%
\end{tabular}%
}%

\newcommand{\BKQUADRIAR}[1]%
{%
\begin{tabular}{c}%
\BKAR{#1}\\%
\rule{0pt}{3pt}\\ \BKAR{#1}\\%
\rule{0pt}{3pt}\\ \BKAR{#1}\\%
\rule{0pt}{3pt}\\ \BKAR{#1}%
\end{tabular}%
}%

\newcommand{\BKQUADRIADJAR}[1]%
{%
\begin{tabular}{c}%
\AR{#1}\\%
\rule{0pt}{3pt}\\ \BKAR{#1}\\%
\rule{0pt}{3pt}\\ \AR{#1}\\%
\rule{0pt}{3pt}\\ \BKAR{#1}%
\end{tabular}%
}%

\newcommand{\BKQUINTIAR}[1]%
{%
\begin{tabular}{c}%
\BKAR{#1}\\%
\rule{0pt}{2pt}\\ \BKAR{#1}\\%
\rule{0pt}{2pt}\\ \BKAR{#1}\\%
\rule{0pt}{2pt}\\ \BKAR{#1}\\%
\rule{0pt}{2pt}\\ \BKAR{#1}%
\end{tabular}%
}%

\newcommand{\BKQUINTIADJAR}[1]%
{%
\begin{tabular}{c}%
\BKAR{#1}\\%
\rule{0pt}{2pt}\\ \AR{#1}\\%
\rule{0pt}{2pt}\\ \BKAR{#1}\\%
\rule{0pt}{2pt}\\ \AR{#1}\\%
\rule{0pt}{2pt}\\ \BKAR{#1}%
\end{tabular}%
}%

\newcommand{\lowername}[2]%
{$\stackrel{\makebox[1pt]{#1}}%
{\begin{picture}(0,0)%
\truex{600}%
\put(0,0){\makebox(0,\value{x})[t]{\makebox[1pt]{$#2$}}}%
\end{picture}}$}%

\newcommand{\hcase}[2]%
{\makebox[0pt]%
{\raisebox{-1pt}[0pt][0pt]{#1{#2}}}}%

\newcommand{\Hcase}[3]%
{\makebox[0pt]
{\raisebox{-1pt}[0pt][0pt]%
{$\stackrel{\makebox[0pt]{$\textstyle{#2}$}}{#1{#3}}$}}}%

\newcommand{\hcasE}[3]%
{\makebox[0pt]%
{\raisebox{-9pt}[0pt][0pt]%
{\lowername{#1{#3}}{#2}}}}%

\newcommand{\hbicase}[2]%
{\makebox[0pt]%
{\raisebox{-2.5pt}[0pt][0pt]{#1{#2}}}}%

\newcommand{\Hbicase}[4]%
{\makebox[0pt]
{\raisebox{-10.5pt}[0pt][0pt]%
{$\stackrel{\makebox[0pt]{$\textstyle{#2}$}}%
{\mbox{\lowername{#1{#4}}{#3}}}$}}}%

\newcommand{\htricase}[2]%
{\makebox[0pt]%
{\raisebox{-13pt}[0pt][0pt]{#1{#2}{}{}{}{}}}}%

\newcommand{\Htricase}[6]%
{\makebox[0pt]%
{\raisebox{-13pt}[0pt][0pt]{#1{#6}{#2}{#3}{#4}{#5}}}}%

\newcommand{\hquadricase}[2]%
{\makebox[0pt]%
{\raisebox{-13pt}[0pt][0pt]{#1{#2}}}}%

\newcommand{\hquinticase}[2]%
{\makebox[0pt]%
{\raisebox{-13pt}[0pt][0pt]{#1{#2}}}}%

\newcommand{\EAR}[1]%
{\begin{picture}(#1,0)%
\put(0,0){\vector(1,0){#1}}%
\end{picture}}%

\newcommand{\EDOTAR}[1]%
{\truex{100}\truey{300}%
\NUMBEROFDOTS=#1%
\divide\NUMBEROFDOTS by \value{y}%
\begin{picture}(#1,0)%
\multiput(0,0)(\value{y},0){\NUMBEROFDOTS}%
{\circle*{\value{x}}}%
\put(#1,0){\vector(1,0){0}}%
\end{picture}}%

\newcommand{\EMONO}[1]%
{\begin{picture}(#1,0)%
\put(0,0){\vector(1,0){#1}}%
\truex{300}\truey{600}%
\put(\value{x},-\value{x}){\line(0,1){\value{y}}}%
\end{picture}}%

\newcommand{\EEPI}[1]%
{\begin{picture}(#1,0)(-#1,0)%
\put(-#1,0){\vector(1,0){#1}}%
\truex{300}\truey{600}\truez{800}%
\put(-\value{z},-\value{x}){\line(0,1){\value{y}}}%
\end{picture}}%

\newcommand{\EBIMO}[1]%
{\begin{picture}(#1,0)(-#1,0)%
\put(-#1,0){\vector(1,0){#1}}%
\truex{300}\truey{600}\truez{800}%
\put(-\value{z},-\value{x}){\line(0,1){\value{y}}}%
\put(-#1,-\value{x}){\hspace{3pt}\line(0,1){\value{y}}}%
\end{picture}}%

\newcommand{\EBIAR}[1]%
{\truex{400}%
\begin{picture}(#1,\value{x})%
\put(0,0){\vector(1,0){#1}}%
\put(0,\value{x}){\vector(1,0){#1}}%
\end{picture}}%

\newcommand{\EEQL}[1]%
{\begin{picture}(#1,0)%
\truex{200}%
\put(0,\value{x}){\line(1,0){#1}}%
\put(0,0){\line(1,0){#1}}%
\end{picture}}%

\newcommand{\EADJAR}[1]%
{\truex{400}%
\begin{picture}(#1,\value{x})%
\put(0,0){\vector(1,0){#1}}%
\put(#1,\value{x}){\vector(-1,0){#1}}%
\end{picture}}%

\newcommand{\ETRIAR}[5]%
{\truex{1200}%
\X=\value{x}%
\multiply\X by 2%
\begin{picture}(#1,\X)%
\put(0,0){\vector(1,0){#1}}%
\put(0,\value{x}){\vector(1,0){#1}}%
\put(0,\X){\vector(1,0){#1}}%
\X=#1%
\divide\X by 2%
\truex{-600}%
\put(\X,\value{x}){\makebox(0,0){$#5$}}%
\truex{600}%
\put(\X,\value{x}){\makebox(0,0){$#4$}}%
\truex{1800}%
\put(\X,\value{x}){\makebox(0,0){$#3$}}%
\truex{3000}%
\put(\X,\value{x}){\makebox(0,0){$#2$}}%
\end{picture}}%

\newcommand{\ETRIADJAR}[5]%
{\truex{1200}%
\X=\value{x}%
\multiply\X by 2%
\begin{picture}(#1,\X)%
\put(#1,0){\vector(-1,0){#1}}%
\put(0,\value{x}){\vector(1,0){#1}}%
\put(#1,\X){\vector(-1,0){#1}}%
\X=#1%
\divide\X by 2%
\truex{-600}%
\put(\X,\value{x}){\makebox(0,0){$#5$}}%
\truex{600}%
\put(\X,\value{x}){\makebox(0,0){$#4$}}%
\truex{1800}%
\put(\X,\value{x}){\makebox(0,0){$#3$}}%
\truex{3000}%
\put(\X,\value{x}){\makebox(0,0){$#2$}}%
\end{picture}}%

\newcommand{\EQUADRIAR}[1]%
{\truex{800}%
\truey{1600}%
\X=\value{x}%
\multiply\X by 3%
\begin{picture}(#1,\X)%
\put(0,0){\vector(1,0){#1}}%
\put(0,\value{x}){\vector(1,0){#1}}%
\put(0,\value{y}){\vector(1,0){#1}}%
\put(0,\X){\vector(1,0){#1}}%
\end{picture}}%

\newcommand{\EQUADRIADJAR}[1]%
{\truex{800}%
\truey{1600}%
\X=\value{x}%
\multiply\X by 3%
\begin{picture}(#1,\X)%
\put(0,0){\vector(1,0){#1}}%
\put(#1,\value{x}){\vector(-1,0){#1}}%
\put(0,\value{y}){\vector(1,0){#1}}%
\put(#1,\X){\vector(-1,0){#1}}%
\end{picture}}%

\newcommand{\EQUINTIAR}[1]%
{\truex{600}%
\truey{1200}%
\truez{1800}%
\X=\value{x}%
\multiply\X by 4%
\begin{picture}(#1,\X)%
\put(0,0){\vector(1,0){#1}}%
\put(0,\value{x}){\vector(1,0){#1}}%
\put(0,\value{y}){\vector(1,0){#1}}%
\put(0,\value{z}){\vector(1,0){#1}}%
\put(0,\X){\vector(1,0){#1}}%
\end{picture}}%

\newcommand{\EQUINTIADJAR}[1]%
{\truex{600}%
\truey{1200}%
\truez{1800}%
\X=\value{x}%
\multiply\X by 4%
\begin{picture}(#1,\X)%
\put(0,0){\vector(1,0){#1}}%
\put(#1,\value{x}){\vector(-1,0){#1}}%
\put(0,\value{y}){\vector(1,0){#1}}%
\put(#1,\value{z}){\vector(-1,0){#1}}%
\put(0,\X){\vector(1,0){#1}}%
\end{picture}}%

\newcommand{\ear}%
{\hspace{\SOURCE\unitlength}%
\hcase{\EAR}{\ARROWLENGTH}}%

\newcommand{\Ear}[1]%
{\hspace{\SOURCE\unitlength}%
\Hcase{\EAR}{#1}{\ARROWLENGTH}}%

\newcommand{\eaR}[1]%
{\hspace{\SOURCE\unitlength}%
\hcasE{\EAR}{#1}{\ARROWLENGTH}}%

\newcommand{\edotar}%
{\hspace{\SOURCE\unitlength}%
\hcase{\EDOTAR}{\ARROWLENGTH}}%

\newcommand{\Edotar}[1]%
{\hspace{\SOURCE\unitlength}%
\Hcase{\EDOTAR}{#1}{\ARROWLENGTH}}%

\newcommand{\edotaR}[1]%
{\hspace{\SOURCE\unitlength}%
\hcasE{\EDOTAR}{#1}{\ARROWLENGTH}}%

\newcommand{\emono}%
{\hspace{\SOURCE\unitlength}%
\hcase{\EMONO}{\ARROWLENGTH}}%

\newcommand{\Emono}[1]%
{\hspace{\SOURCE\unitlength}%
\Hcase{\EMONO}{#1}{\ARROWLENGTH}}%

\newcommand{\emonO}[1]%
{\hspace{\SOURCE\unitlength}%
\hcasE{\EMONO}{#1}{\ARROWLENGTH}}%

\newcommand{\eepi}%
{\hspace{\SOURCE\unitlength}%
\hcase{\EEPI}{\ARROWLENGTH}}%

\newcommand{\Eepi}[1]%
{\hspace{\SOURCE\unitlength}%
\Hcase{\EEPI}{#1}{\ARROWLENGTH}}%

\newcommand{\eepI}[1]%
{\hspace{\SOURCE\unitlength}%
\hcasE{\EEPI}{#1}{\ARROWLENGTH}}%

\newcommand{\ebimo}%
{\hspace{\SOURCE\unitlength}%
\hcase{\EBIMO}{\ARROWLENGTH}}%

\newcommand{\Ebimo}[1]%
{\hspace{\SOURCE\unitlength}%
\Hcase{\EBIMO}{#1}{\ARROWLENGTH}}%

\newcommand{\ebimO}[1]%
{\hspace{\SOURCE\unitlength}%
\hcasE{\EBIMO}{#1}{\ARROWLENGTH}}%

\newcommand{\eiso}%
{\hspace{\SOURCE\unitlength}%
\Hcase{\EAR}{\cong}{\ARROWLENGTH}}%

\newcommand{\Eiso}[1]%
{\hspace{\SOURCE\unitlength}%
\Hcase{\EAR}{\cong#1}{\ARROWLENGTH}}%

\newcommand{\eisO}[1]%
{\hspace{\SOURCE\unitlength}%
\hcasE{\EAR}{\cong#1}{\ARROWLENGTH}}%

\newcommand{\ebiar}%
{\hspace{\SOURCE\unitlength}%
\hbicase{\EBIAR}{\ARROWLENGTH}}%

\newcommand{\Ebiar}[2]%
{\hspace{\SOURCE\unitlength}%
\Hbicase{\EBIAR}{#1}{#2}{\ARROWLENGTH}}%

\newcommand{\eeql}%
{\hspace{\SOURCE\unitlength}%
\hbicase{\EEQL}{\ARROWLENGTH}}%

\newcommand{\eadjar}%
{\hspace{\SOURCE\unitlength}%
\hbicase{\EADJAR}{\ARROWLENGTH}}%

\newcommand{\Eadjar}[2]%
{\hspace{\SOURCE\unitlength}%
\Hbicase{\EADJAR}{#1}{#2}{\ARROWLENGTH}}%

\newcommand{\etriar}%
{\hspace{\SOURCE\unitlength}%
\htricase{\ETRIAR}{\ARROWLENGTH}}%

\newcommand{\Etriar}[3]%
{\hspace{\SOURCE\unitlength}%
\Htricase{\ETRIAR}{#1}{#2}{#3}{}{\ARROWLENGTH}}%

\newcommand{\etriaR}[3]%
{\hspace{\SOURCE\unitlength}%
\Htricase{\ETRIAR}{}{#1}{#2}{#3}{\ARROWLENGTH}}%

\newcommand{\etriadjar}%
{\hspace{\SOURCE\unitlength}%
\htricase{\ETRIADJAR}{\ARROWLENGTH}}%

\newcommand{\Etriadjar}[3]%
{\hspace{\SOURCE\unitlength}%
\Htricase{\ETRIADJAR}{#1}{#2}{#3}{}{\ARROWLENGTH}}%

\newcommand{\etriadjaR}[3]%
{\hspace{\SOURCE\unitlength}%
\Htricase{\ETRIADJAR}{}{#1}{#2}{#3}{\ARROWLENGTH}}%

\newcommand{\equadriar}%
{\hspace{\SOURCE\unitlength}%
\hquadricase{\EQUADRIAR}{\ARROWLENGTH}}%

\newcommand{\equadriadjar}%
{\hspace{\SOURCE\unitlength}%
\hquadricase{\EQUADRIADJAR}{\ARROWLENGTH}}%

\newcommand{\equintiar}%
{\hspace{\SOURCE\unitlength}%
\hquinticase{\EQUINTIAR}{\ARROWLENGTH}}%

\newcommand{\equintiadjar}%
{\hspace{\SOURCE\unitlength}%
\hquinticase{\EQUINTIADJAR}{\ARROWLENGTH}}%

\newcommand{\WAR}[1]%
{\begin{picture}(#1,0)%
\put(#1,0){\vector(-1,0){#1}}%
\end{picture}}%

\newcommand{\WDOTAR}[1]%
{\truex{100}\truey{300}%
\NUMBEROFDOTS=#1%
\divide\NUMBEROFDOTS by \value{y}%
\begin{picture}(#1,0)%
\multiput(#1,0)(-\value{y},0){\NUMBEROFDOTS}%
{\circle*{\value{x}}}%
\put(0,0){\vector(-1,0){0}}%
\end{picture}}%

\newcommand{\WMONO}[1]%
{\begin{picture}(#1,0)(-#1,0)%
\put(0,0){\vector(-1,0){#1}}%
\truex{300}\truey{600}%
\put(-\value{x},-\value{x}){\line(0,1){\value{y}}}%
\end{picture}}%

\newcommand{\WEPI}[1]%
{\begin{picture}(#1,0)%
\put(#1,0){\vector(-1,0){#1}}%
\truex{300}\truey{600}\truez{800}%
\put(\value{z},-\value{x}){\line(0,1){\value{y}}}%
\end{picture}}%

\newcommand{\WBIMO}[1]%
{\begin{picture}(#1,0)%
\put(#1,0){\vector(-1,0){#1}}%
\truex{300}\truey{600}\truez{800}%
\put(\value{z},-\value{x}){\line(0,1){\value{y}}}%
\put(#1,-\value{x}){\hspace{-3pt}\line(0,1){\value{y}}}%
\end{picture}}%

\newcommand{\WBIAR}[1]%
{\truex{400}%
\begin{picture}(#1,\value{x})%
\put(#1,0){\vector(-1,0){#1}}%
\put(#1,\value{x}){\vector(-1,0){#1}}%
\end{picture}}%

\newcommand{\WADJAR}[1]%
{\truex{400}%
\begin{picture}(#1,\value{x})%
\put(0,\value{x}){\vector(1,0){#1}}%
\put(#1,0){\vector(-1,0){#1}}%
\end{picture}}%

\newcommand{\WTRIAR}[5]%
{\truex{1200}%
\X=\value{x}%
\multiply\X by 2%
\begin{picture}(#1,\X)%
\put(#1,0){\vector(-1,0){#1}}%
\put(#1,\value{x}){\vector(-1,0){#1}}%
\put(#1,\X){\vector(-1,0){#1}}%
\X=#1%
\divide\X by 2%
\truex{-600}%
\put(\X,\value{x}){\makebox(0,0){$#5$}}%
\truex{600}%
\put(\X,\value{x}){\makebox(0,0){$#4$}}%
\truex{1800}%
\put(\X,\value{x}){\makebox(0,0){$#3$}}%
\truex{3000}%
\put(\X,\value{x}){\makebox(0,0){$#2$}}%
\end{picture}}%

\newcommand{\WTRIADJAR}[5]%
{\truex{1200}%
\X=\value{x}%
\multiply\X by 2%
\begin{picture}(#1,\X)%
\put(0,0){\vector(1,0){#1}}%
\put(#1,\value{x}){\vector(-1,0){#1}}%
\put(0,\X){\vector(1,0){#1}}%
\X=#1%
\divide\X by 2%
\truex{-600}%
\put(\X,\value{x}){\makebox(0,0){$#5$}}%
\truex{600}%
\put(\X,\value{x}){\makebox(0,0){$#4$}}%
\truex{1800}%
\put(\X,\value{x}){\makebox(0,0){$#3$}}%
\truex{3000}%
\put(\X,\value{x}){\makebox(0,0){$#2$}}%
\end{picture}}%

\newcommand{\WQUADRIAR}[1]%
{\truex{800}%
\truey{1600}%
\X=\value{x}%
\multiply\X by 3%
\begin{picture}(#1,\X)%
\put(#1,0){\vector(-1,0){#1}}%
\put(#1,\value{x}){\vector(-1,0){#1}}%
\put(#1,\value{y}){\vector(-1,0){#1}}%
\put(#1,\X){\vector(-1,0){#1}}%
\end{picture}}%

\newcommand{\WQUADRIADJAR}[1]%
{\truex{800}%
\truey{1600}%
\X=\value{x}%
\multiply\X by 3%
\begin{picture}(#1,\X)%
\put(#1,0){\vector(-1,0){#1}}%
\put(0,\value{x}){\vector(1,0){#1}}%
\put(#1,\value{y}){\vector(-1,0){#1}}%
\put(0,\X){\vector(1,0){#1}}%
\end{picture}}%

\newcommand{\WQUINTIAR}[1]%
{\truex{600}%
\truey{1200}%
\truez{1800}%
\X=\value{x}%
\multiply\X by 4%
\begin{picture}(#1,\X)%
\put(#1,0){\vector(-1,0){#1}}%
\put(#1,\value{x}){\vector(-1,0){#1}}%
\put(#1,\value{y}){\vector(-1,0){#1}}%
\put(#1,\value{z}){\vector(-1,0){#1}}%
\put(#1,\X){\vector(-1,0){#1}}%
\end{picture}}%

\newcommand{\WQUINTIADJAR}[1]%
{\truex{600}%
\truey{1200}%
\truez{1800}%
\X=\value{x}%
\multiply\X by 4%
\begin{picture}(#1,\X)%
\put(#1,0){\vector(-1,0){#1}}%
\put(0,\value{x}){\vector(1,0){#1}}%
\put(#1,\value{y}){\vector(-1,0){#1}}%
\put(0,\value{z}){\vector(1,0){#1}}%
\put(#1,\X){\vector(-1,0){#1}}%
\end{picture}}%

\newcommand{\war}%
{\hspace{\SOURCE\unitlength}%
\hcase{\WAR}{\ARROWLENGTH}}%

\newcommand{\War}[1]%
{\hspace{\SOURCE\unitlength}%
\Hcase{\WAR}{#1}{\ARROWLENGTH}}%

\newcommand{\waR}[1]%
{\hspace{\SOURCE\unitlength}%
\hcasE{\WAR}{#1}{\ARROWLENGTH}}%

\newcommand{\wdotar}%
{\hspace{\SOURCE\unitlength}%
\hcase{\WDOTAR}{\ARROWLENGTH}}%

\newcommand{\Wdotar}[1]%
{\hspace{\SOURCE\unitlength}%
\Hcase{\WDOTAR}{#1}{\ARROWLENGTH}}%

\newcommand{\wdotaR}[1]%
{\hspace{\SOURCE\unitlength}%
\hcasE{\WDOTAR}{#1}{\ARROWLENGTH}}%

\newcommand{\wmono}%
{\hspace{\SOURCE\unitlength}%
\hcase{\WMONO}{\ARROWLENGTH}}%

\newcommand{\Wmono}[1]%
{\hspace{\SOURCE\unitlength}%
\Hcase{\WMONO}{#1}{\ARROWLENGTH}}%

\newcommand{\wmonO}[1]%
{\hspace{\SOURCE\unitlength}%
\hcasE{\WMONO}{#1}{\ARROWLENGTH}}%

\newcommand{\wepi}%
{\hspace{\SOURCE\unitlength}%
\hcase{\WEPI}{\ARROWLENGTH}}%

\newcommand{\Wepi}[1]%
{\hspace{\SOURCE\unitlength}%
\Hcase{\WEPI}{#1}{\ARROWLENGTH}}%

\newcommand{\wepI}[1]%
{\hspace{\SOURCE\unitlength}%
\hcasE{\WEPI}{#1}{\ARROWLENGTH}}%

\newcommand{\wbimo}%
{\hspace{\SOURCE\unitlength}%
\hcase{\WBIMO}{\ARROWLENGTH}}%

\newcommand{\Wbimo}[1]%
{\hspace{\SOURCE\unitlength}%
\Hcase{\WBIMO}{#1}{\ARROWLENGTH}}%

\newcommand{\wbimO}[1]%
{\hspace{\SOURCE\unitlength}%
\hcasE{\WBIMO}{#1}{\ARROWLENGTH}}%

\newcommand{\wiso}%
{\hspace{\SOURCE\unitlength}%
\Hcase{\WAR}{\cong}{\ARROWLENGTH}}%

\newcommand{\Wiso}[1]%
{\hspace{\SOURCE\unitlength}%
\Hcase{\WAR}{\cong#1}{\ARROWLENGTH}}%

\newcommand{\wisO}[1]%
{\hspace{\SOURCE\unitlength}%
\hcasE{\WAR}{#1}{\ARROWLENGTH}}%

\newcommand{\wbiar}%
{\hspace{\SOURCE\unitlength}%
\hbicase{\WBIAR}{\ARROWLENGTH}}%

\newcommand{\Wbiar}[2]%
{\hspace{\SOURCE\unitlength}%
\Hbicase{\WBIAR}{#1}{#2}{\ARROWLENGTH}}%

\newcommand{\weql}%
{\hspace{\SOURCE\unitlength}%
\hbicase{\EEQL}{\ARROWLENGTH}}%

\newcommand{\wadjar}%
{\hspace{\SOURCE\unitlength}%
\hbicase{\WADJAR}{\ARROWLENGTH}}%

\newcommand{\Wadjar}[2]%
{\hspace{\SOURCE\unitlength}%
\Hbicase{\WADJAR}{#1}{#2}{\ARROWLENGTH}}%

\newcommand{\wtriar}%
{\hspace{\SOURCE\unitlength}%
\htricase{\WTRIAR}{\ARROWLENGTH}}%

\newcommand{\Wtriar}[3]%
{\hspace{\SOURCE\unitlength}%
\Htricase{\WTRIAR}{#1}{#2}{#3}{}{\ARROWLENGTH}}%

\newcommand{\wtriaR}[3]%
{\hspace{\SOURCE\unitlength}%
\Htricase{\WTRIAR}{}{#1}{#2}{#3}{\ARROWLENGTH}}%

\newcommand{\wtriadjar}%
{\hspace{\SOURCE\unitlength}%
\htricase{\WTRIADJAR}{\ARROWLENGTH}}%

\newcommand{\Wtriadjar}[3]%
{\hspace{\SOURCE\unitlength}%
\Htricase{\WTRIADJAR}{#1}{#2}{#3}{}{\ARROWLENGTH}}%

\newcommand{\wtriadjaR}[3]%
{\hspace{\SOURCE\unitlength}%
\Htricase{\WTRIADJAR}{}{#1}{#2}{#3}{\ARROWLENGTH}}%

\newcommand{\wquadriar}%
{\hspace{\SOURCE\unitlength}%
\hquadricase{\WQUADRIAR}{\ARROWLENGTH}}%

\newcommand{\wquadriadjar}%
{\hspace{\SOURCE\unitlength}%
\hquadricase{\WQUADRIADJAR}{\ARROWLENGTH}}%

\newcommand{\wquintiar}%
{\hspace{\SOURCE\unitlength}%
\hquinticase{\WQUINTIAR}{\ARROWLENGTH}}%

\newcommand{\wquintiadjar}%
{\hspace{\SOURCE\unitlength}%
\hquinticase{\WQUINTIADJAR}{\ARROWLENGTH}}%

\newcommand{\Vcase}[3]{\makebox[0pt]%
{\makebox[0pt][r]{\raisebox{0pt}[0pt][0pt]{${#2}\hspace{2pt}$}}}#1{#3}}%

\newcommand{\vcasE}[3]{\makebox[0pt]%
{#1{#3}\makebox[0pt][l]{\raisebox{0pt}[0pt][0pt]{\hspace{2pt}$#2$}}}}%

\newcommand{\Vbicase}[4]{\makebox[0pt]%
{\makebox[0pt][r]{\raisebox{0pt}[0pt][0pt]{$#2$\hspace{4pt}}}#1{#4}%
\makebox[0pt][l]{\raisebox{0pt}[0pt][0pt]{\hspace{5pt}$#3$}}}}%

\newcommand{\vtricase}[2]%
{\makebox[0pt]{#1{}{}{}{}{#2}}}%

\newcommand{\Vtricase}[6]%
{\makebox[0pt]{#1{#2}{#3}{#4}{#5}{#6}}}%

\newcommand{\SAR}[1]%
{\begin{picture}(0,0)%
\put(0,0){\makebox(0,0)%
{\begin{picture}(0,#1)%
\put(0,#1){\vector(0,-1){#1}}%
\end{picture}}}\end{picture}}%

\newcommand{\SDOTAR}[1]%
{\truex{100}\truey{300}%
\NUMBEROFDOTS=#1%
\divide\NUMBEROFDOTS by \value{y}%
\begin{picture}(0,0)%
\put(0,0){\makebox(0,0)%
{\begin{picture}(0,#1)%
\multiput(0,#1)(0,-\value{y}){\NUMBEROFDOTS}%
{\circle*{\value{x}}}%
\put(0,0){\vector(0,-1){0}}%
\end{picture}}}\end{picture}}%

\newcommand{\SMONO}[1]%
{\begin{picture}(0,0)%
\put(0,0){\makebox(0,0)%
{\begin{picture}(0,#1)%
\put(0,#1){\vector(0,-1){#1}}%
\truex{300}\truey{600}%
\put(0,#1){\begin{picture}(0,0)%
\put(-\value{x},-\value{x}){\line(1,0){\value{y}}}\end{picture}}%
\end{picture}}}\end{picture}}%

\newcommand{\SEPI}[1]%
{\begin{picture}(0,0)%
\put(0,0){\makebox(0,0)%
{\begin{picture}(0,#1)%
\put(0,#1){\vector(0,-1){#1}}%
\truex{300}\truey{600}\truez{800}%
\put(-\value{x},\value{z}){\line(1,0){\value{y}}}%
\end{picture}}}\end{picture}}%

\newcommand{\SBIMO}[1]%
{\begin{picture}(0,0)%
\put(0,0){\makebox(0,0)%
{\begin{picture}(0,#1)%
\put(0,#1){\vector(0,-1){#1}}%
\truex{300}\truey{600}\truez{800}%
\put(0,#1){\begin{picture}(0,0)%
\put(-\value{x},-\value{x}){\line(1,0){\value{y}}}\end{picture}}%
\put(-\value{x},\value{z}){\line(1,0){\value{y}}}%
\end{picture}}}\end{picture}}%

\newcommand{\SBIAR}[1]%
{\begin{picture}(0,0)%
\truex{200}%
\put(0,0){\makebox(0,0)%
{\begin{picture}(0,#1)\put(-\value{x},#1){\vector(0,-1){#1}}%
\put(\value{x},#1){\vector(0,-1){#1}}%
\end{picture}}}\end{picture}}%

\newcommand{\SEQL}[1]%
{\begin{picture}(0,0)%
\truex{100}%
\put(0,0){\makebox(0,0)%
{\begin{picture}(0,#1)\put(-\value{x},#1){\line(0,-1){#1}}%
\put(\value{x},#1){\line(0,-1){#1}}%
\end{picture}}}\end{picture}}%

\newcommand{\STRIAR}[5]%
{\truex{1200}%
\truey{1800}%
\truez{600}%
\begin{picture}(0,0)%
\put(0,0){\makebox(0,0)%
{\begin{picture}(0,#5)%
\X=#5\divide\X by 2%
\put(-\value{x},#5){\vector(0,-1){#5}}%
\put(0,#5){\vector(0,-1){#5}}%
\put(\value{x},#5){\vector(0,-1){#5}}%
\put(-\value{y},\X){\makebox(0,0){#1}}%
\put(\value{y},\X){\makebox(0,0){#4}}%
\put(-\value{z},\X){\makebox(0,0){#2}}%
\put(\value{z},\X){\makebox(0,0){#3}}%
\end{picture}}}\end{picture}}%

\newcommand{\STRIADJAR}[5]%
{\truex{1200}%
\truey{1800}%
\truez{600}%
\begin{picture}(0,0)%
\put(0,0){\makebox(0,0)%
{\begin{picture}(0,#5)%
\X=#5\divide\X by 2%
\put(-\value{x},0){\vector(0,1){#5}}%
\put(0,#5){\vector(0,-1){#5}}%
\put(\value{x},0){\vector(0,1){#5}}%
\put(-\value{y},\X){\makebox(0,0){#1}}%
\put(\value{y},\X){\makebox(0,0){#4}}%
\put(-\value{z},\X){\makebox(0,0){#2}}%
\put(\value{z},\X){\makebox(0,0){#3}}%
\end{picture}}}\end{picture}}%

\newcommand{\SQUADRIAR}[1]%
{\truex{400}%
\truey{1200}%
\begin{picture}(0,0)%
\put(0,0){\makebox(0,0)%
{\begin{picture}(0,#1)%
\put(-\value{y},#1){\vector(0,-1){#1}}%
\put(-\value{x},#1){\vector(0,-1){#1}}%
\put(\value{x},#1){\vector(0,-1){#1}}%
\put(\value{y},#1){\vector(0,-1){#1}}%
\end{picture}}}\end{picture}}%

\newcommand{\SQUADRIADJAR}[1]%
{\truex{400}%
\truey{1200}%
\begin{picture}(0,0)%
\put(0,0){\makebox(0,0)%
{\begin{picture}(0,#1)%
\put(-\value{y},#1){\vector(0,-1){#1}}%
\put(-\value{x},0){\vector(0,1){#1}}%
\put(\value{x},#1){\vector(0,-1){#1}}%
\put(\value{y},0){\vector(0,1){#1}}%
\end{picture}}}\end{picture}}%

\newcommand{\SQUINTIAR}[1]%
{\truex{600}%
\truey{1200}%
\begin{picture}(0,0)%
\put(0,0){\makebox(0,0)%
{\begin{picture}(0,#1)%
\put(-\value{y},#1){\vector(0,-1){#1}}%
\put(-\value{x},#1){\vector(0,-1){#1}}%
\put(0,#1){\vector(0,-1){#1}}%
\put(\value{x},#1){\vector(0,-1){#1}}%
\put(\value{y},#1){\vector(0,-1){#1}}%
\end{picture}}}\end{picture}}%

\newcommand{\SQUINTIADJAR}[1]%
{\truex{600}%
\truey{1200}%
\begin{picture}(0,0)%
\put(0,0){\makebox(0,0)%
{\begin{picture}(0,#1)%
\put(-\value{y},#1){\vector(0,-1){#1}}%
\put(-\value{x},0){\vector(0,1){#1}}%
\put(0,#1){\vector(0,-1){#1}}%
\put(\value{x},0){\vector(0,1){#1}}%
\put(\value{y},#1){\vector(0,-1){#1}}%
\end{picture}}}\end{picture}}%

\newcommand{\Sarv}[2]{\Vcase{\SAR}{#1}{#200}}%

\newcommand{\Sar}[1]{\Sarv{#1}{50}}%

\newcommand{\saRv}[2]{\vcasE{\SAR}{#1}{#200}}%

\newcommand{\saR}[1]{\saRv{#1}{50}}%

\newcommand{\Sisov}[2]%
{\Vbicase{\SAR}{#1\hspace{-2pt}}{\hspace{-2pt}\cong}{#200}}%

\newcommand{\NAR}[1]%
{\begin{picture}(0,0)%
\put(0,0){\makebox(0,0)%
{\begin{picture}(0,#1)\put(0,0){\vector(0,1){#1}}%
\end{picture}}}\end{picture}}%

\newcommand{\NDOTAR}[1]%
{\truex{100}\truey{300}%
\NUMBEROFDOTS=#1%
\divide\NUMBEROFDOTS by \value{y}%
\begin{picture}(0,0)%
\put(0,0){\makebox(0,0)%
{\begin{picture}(0,#1)%
\multiput(0,0)(0,\value{y}){\NUMBEROFDOTS}%
{\circle*{\value{x}}}%
\put(0,#1){\vector(0,1){0}}%
\end{picture}}}\end{picture}}%

\newcommand{\NMONO}[1]%
{\begin{picture}(0,0)%
\put(0,0){\makebox(0,0)%
{\begin{picture}(0,#1)%
\put(0,0){\vector(0,1){#1}}%
\truex{300}\truey{600}%
\put(-\value{x},\value{x}){\line(1,0){\value{y}}}%
\end{picture}}}%
\end{picture}}%

\newcommand{\NEPI}[1]%
{\begin{picture}(0,0)%
\put(0,0){\makebox(0,0)%
{\begin{picture}(0,#1)%
\put(0,0){\vector(0,1){#1}}%
\truex{300}\truey{600}\truez{800}%
\put(0,#1){\begin{picture}(0,0)%
\put(-\value{x},-\value{z}){\line(1,0){\value{y}}}\end{picture}}%
\end{picture}}}\end{picture}}%

\newcommand{\NBIMO}[1]%
{\begin{picture}(0,0)%
\put(0,0){\makebox(0,0)%
{\begin{picture}(0,#1)%
\put(0,0){\vector(0,1){#1}}%
\truex{300}\truey{600}\truez{800}%
\put(-\value{x},\value{x}){\line(1,0){\value{y}}}%
\put(0,#1){\begin{picture}(0,0)%
\put(-\value{x},-\value{z}){\line(1,0){\value{y}}}\end{picture}}%
\end{picture}}}\end{picture}}%

\newcommand{\NBIAR}[1]%
{\begin{picture}(0,0)%
\truex{200}%
\put(0,0){\makebox(0,0)%
{\begin{picture}(0,#1)\put(-\value{x},0){\vector(0,1){#1}}%
\put(\value{x},0){\vector(0,1){#1}}%
\end{picture}}}\end{picture}}%

\newcommand{\NTRIAR}[5]%
{\truex{1200}%
\truey{1800}%
\truez{600}%
\begin{picture}(0,0)%
\put(0,0){\makebox(0,0)%
{\begin{picture}(0,#5)%
\X=#5\divide\X by 2%
\put(-\value{x},0){\vector(0,1){#5}}%
\put(0,0){\vector(0,1){#5}}%
\put(\value{x},0){\vector(0,1){#5}}%
\put(-\value{y},\X){\makebox(0,0){#1}}%
\put(\value{y},\X){\makebox(0,0){#4}}%
\put(-\value{z},\X){\makebox(0,0){#2}}%
\put(\value{z},\X){\makebox(0,0){#3}}%
\end{picture}}}\end{picture}}%

\newcommand{\NTRIADJAR}[5]%
{\truex{1200}%
\truey{1800}%
\truez{600}%
\begin{picture}(0,0)%
\put(0,0){\makebox(0,0)%
{\begin{picture}(0,#5)%
\X=#5\divide\X by 2%
\put(-\value{x},#5){\vector(0,-1){#5}}%
\put(0,0){\vector(0,1){#5}}%
\put(\value{x},#5){\vector(0,-1){#5}}%
\put(-\value{y},\X){\makebox(0,0){#1}}%
\put(\value{y},\X){\makebox(0,0){#4}}%
\put(-\value{z},\X){\makebox(0,0){#2}}%
\put(\value{z},\X){\makebox(0,0){#3}}%
\end{picture}}}\end{picture}}%

\newcommand{\NQUADRIAR}[1]%
{\truex{400}%
\truey{1200}%
\begin{picture}(0,0)%
\put(0,0){\makebox(0,0)%
{\begin{picture}(0,#1)%
\put(-\value{y},0){\vector(0,1){#1}}%
\put(-\value{x},0){\vector(0,1){#1}}%
\put(\value{x},0){\vector(0,1){#1}}%
\put(\value{y},0){\vector(0,1){#1}}%
\end{picture}}}\end{picture}}%

\newcommand{\NQUADRIADJAR}[1]%
{\truex{400}%
\truey{1200}%
\begin{picture}(0,0)%
\put(0,0){\makebox(0,0)%
{\begin{picture}(0,#1)%
\put(-\value{y},0){\vector(0,1){#1}}%
\put(-\value{x},#1){\vector(0,-1){#1}}%
\put(\value{x},0){\vector(0,1){#1}}%
\put(\value{y},#1){\vector(0,-1){#1}}%
\end{picture}}}\end{picture}}%

\newcommand{\NQUINTIAR}[1]%
{\truex{600}%
\truey{1200}%
\begin{picture}(0,0)%
\put(0,0){\makebox(0,0)%
{\begin{picture}(0,#1)%
\put(-\value{y},0){\vector(0,1){#1}}%
\put(-\value{x},0){\vector(0,1){#1}}%
\put(0,0){\vector(0,1){#1}}%
\put(\value{x},0){\vector(0,1){#1}}%
\put(\value{y},0){\vector(0,1){#1}}%
\end{picture}}}\end{picture}}%

\newcommand{\NQUINTIADJAR}[1]%
{\truex{600}%
\truey{1200}%
\begin{picture}(0,0)%
\put(0,0){\makebox(0,0)%
{\begin{picture}(0,#1)%
\put(-\value{y},0){\vector(0,1){#1}}%
\put(-\value{x},#1){\vector(0,-1){#1}}%
\put(0,0){\vector(0,1){#1}}%
\put(\value{x},#1){\vector(0,-1){#1}}%
\put(\value{y},0){\vector(0,1){#1}}%
\end{picture}}}\end{picture}}%

\newcommand{\Nisov}[2]%
{\Vbicase{\NAR}{#1\hspace{-2pt}}{\hspace{-2pt}\cong}{#200}}%

\newcommand{\NEDOTAR}%
{\truex{100}\truey{212}%
\NUMBEROFDOTS=5800%
\divide\NUMBEROFDOTS by \value{y}%
\begin{picture}(0,0)%
\multiput(-2900,-2900)(\value{y},\value{y}){\NUMBEROFDOTS}%
{\circle*{\value{x}}}%
\put(2900,2900){\vector(1,1){0}}%
\end{picture}}%

\newcommand{\SWAR}{\begin{picture}(0,0)%
\put(2900,2900){\vector(-1,-1){5800}}%
\end{picture}}%

\newcommand{\SWDOTAR}%
{\truex{100}\truey{212}%
\NUMBEROFDOTS=5800%
\divide\NUMBEROFDOTS by \value{y}%
\begin{picture}(0,0)%
\multiput(2900,2900)(-\value{y},-\value{y}){\NUMBEROFDOTS}%
{\circle*{\value{x}}}%
\put(-2900,-2900){\vector(-1,-1){0}}%
\end{picture}}%

\newcommand{\SEDOTAR}%
{\truex{100}\truey{212}%
\NUMBEROFDOTS=5800%
\divide\NUMBEROFDOTS by \value{y}%
\begin{picture}(0,0)%
\multiput(-2900,2900)(\value{y},-\value{y}){\NUMBEROFDOTS}%
{\circle*{\value{x}}}%
\put(2900,-2900){\vector(1,-1){0}}%
\end{picture}}%

\newcommand{\NWDOTAR}%
{\truex{100}\truey{212}%
\NUMBEROFDOTS=5800%
\divide\NUMBEROFDOTS by \value{y}%
\begin{picture}(0,0)%
\multiput(2900,-2900)(-\value{y},\value{y}){\NUMBEROFDOTS}%
{\circle*{\value{x}}}%
\put(-2900,2900){\vector(-1,1){0}}%
\end{picture}}%

\newcommand{\ENEAR}[2]%
{\makebox[0pt]{\begin{picture}(0,0)%
\put(0,-150){\makebox(0,0){\begin{picture}(0,0)%
\put(-6600,-3300){\vector(2,1){13200}}%
\truex{200}\truey{800}\truez{600}%
\put(-\value{x},\value{x}){\makebox(0,\value{z})[r]{${#1}$}}%
\put(\value{x},-\value{y}){\makebox(0,\value{z})[l]{${#2}$}}%
\end{picture}}}\end{picture}}}%

\newcommand{\ESEAR}[2]%
{\makebox[0pt]{\begin{picture}(0,0)%
\put(0,-150){\makebox(0,0){\begin{picture}(0,0)%
\put(-6600,3300){\vector(2,-1){13200}}%
\truex{200}\truey{800}\truez{600}%
\put(\value{x},\value{x}){\makebox(0,\value{z})[l]{${#1}$}}%
\put(-\value{x},-\value{y}){\makebox(0,\value{z})[r]{${#2}$}}%
\end{picture}}}\end{picture}}}%

\newcommand{\WNWAR}[2]%
{\makebox[0pt]{\begin{picture}(0,0)%
\put(0,-150){\makebox(0,0){\begin{picture}(0,0)%
\put(6600,-3300){\vector(-2,1){13200}}%
\truex{200}\truey{800}\truez{600}%
\put(\value{x},\value{x}){\makebox(0,\value{z})[l]{${#1}$}}%
\put(-\value{x},-\value{y}){\makebox(0,\value{z})[r]{${#2}$}}%
\end{picture}}}\end{picture}}}%

\newcommand{\WSWAR}[2]%
{\makebox[0pt]{\begin{picture}(0,0)%
\put(0,-150){\makebox(0,0){\begin{picture}(0,0)%
\put(6600,3300){\vector(-2,-1){13200}}%
\truex{200}\truey{800}\truez{600}%
\put(-\value{x},\value{x}){\makebox(0,\value{z})[r]{${#1}$}}%
\put(\value{x},-\value{y}){\makebox(0,\value{z})[l]{${#2}$}}%
\end{picture}}}\end{picture}}}%

\newcommand{\NNEAR}[2]%
{\raisebox{-1pt}[0pt][0pt]{\begin{picture}(0,0)%
\put(0,0){\makebox(0,0){\begin{picture}(0,0)%
\put(-3300,-6600){\vector(1,2){6600}}%
\truex{100}\truez{600}%
\put(-\value{x},\value{x}){\makebox(0,\value{z})[r]{${#1}$}}%
\put(\value{x},-\value{z}){\makebox(0,\value{z})[l]{${#2}$}}%
\end{picture}}}\end{picture}}}%

\newcommand{\SSWAR}[2]%
{\raisebox{-1pt}[0pt][0pt]{\begin{picture}(0,0)%
\put(0,0){\makebox(0,0){\begin{picture}(0,0)%
\put(3300,6600){\vector(-1,-2){6600}}%
\truex{100}\truez{600}%
\put(-\value{x},\value{x}){\makebox(0,\value{z})[r]{${#1}$}}%
\put(\value{x},-\value{z}){\makebox(0,\value{z})[l]{${#2}$}}%
\end{picture}}}\end{picture}}}%

\newcommand{\SSEAR}[2]%
{\raisebox{-1pt}[0pt][0pt]{\begin{picture}(0,0)%
\put(0,0){\makebox(0,0){\begin{picture}(0,0)%
\put(-3300,6600){\vector(1,-2){6600}}%
\truex{200}\truez{600}%
\put(\value{x},\value{x}){\makebox(0,\value{z})[l]{${#1}$}}%
\put(-\value{x},-\value{z}){\makebox(0,\value{z})[r]{${#2}$}}%
\end{picture}}}\end{picture}}}%

\newcommand{\NNWAR}[2]%
{\raisebox{-1pt}[0pt][0pt]{\begin{picture}(0,0)%
\put(0,0){\makebox(0,0){\begin{picture}(0,0)%
\put(3300,-6600){\vector(-1,2){6600}}%
\truex{200}\truez{600}%
\put(\value{x},\value{x}){\makebox(0,\value{z})[l]{${#1}$}}%
\put(-\value{x},-\value{z}){\makebox(0,\value{z})[r]{${#2}$}}%
\end{picture}}}\end{picture}}}%

\newcommand{\Necurve}[2]%
{\begin{picture}(0,0)%
\truex{1300}\truey{2000}\truez{200}%
\put(0,\value{x}){\oval(#200,\value{y})[t]}%
\put(0,\value{x}){\makebox(0,0){\begin{picture}(#200,0)%
\put(#200,0){\vector(0,-1){\value{z}}}%
\put(0,0){\line(0,-1){\value{z}}}\end{picture}}}%
\truex{2500}%
\put(0,\value{x}){\makebox(0,0)[b]{${#1}$}}%
\end{picture}}%

\newcommand{\Nwcurve}[2]%
{\begin{picture}(0,0)%
\truex{1300}\truey{2000}\truez{200}%
\put(0,\value{x}){\oval(#200,\value{y})[t]}%
\put(0,\value{x}){\makebox(0,0){\begin{picture}(#200,0)%
\put(#200,0){\line(0,-1){\value{z}}}%
\put(0,0){\vector(0,-1){\value{z}}}\end{picture}}}%
\truex{2500}%
\put(0,\value{x}){\makebox(0,0)[b]{${#1}$}}%
\end{picture}}%

\newcommand{\Securve}[2]%
{\begin{picture}(0,0)%
\truex{1300}\truey{2000}\truez{200}%
\put(0,-\value{x}){\oval(#200,\value{y})[b]}%
\put(0,-\value{x}){\makebox(0,0){\begin{picture}(#200,0)%
\put(#200,0){\vector(0,1){\value{z}}}%
\put(0,0){\line(0,1){\value{z}}}\end{picture}}}%
\truex{2500}%
\put(0,-\value{x}){\makebox(0,0)[t]{${#1}$}}%
\end{picture}}%

\newcommand{\Swcurve}[2]%
{\begin{picture}(0,0)%
\truex{1300}\truey{2000}\truez{200}%
\put(0,-\value{x}){\oval(#200,\value{y})[b]}%
\put(0,-\value{x}){\makebox(0,0){\begin{picture}(#200,0)%
\put(#200,0){\line(0,1){\value{z}}}%
\put(0,0){\vector(0,1){\value{z}}}\end{picture}}}%
\truex{2500}%
\put(0,-\value{x}){\makebox(0,0)[t]{${#1}$}}%
\end{picture}}%

\newcommand{\Escurve}[2]%
{\begin{picture}(0,0)%
\truex{1400}\truey{2000}\truez{200}%
\put(\value{x},0){\oval(\value{y},#200)[r]}%
\put(\value{x},0){\makebox(0,0){\begin{picture}(0,#200)%
\put(0,0){\vector(-1,0){\value{z}}}%
\put(0,#200){\line(-1,0){\value{z}}}\end{picture}}}%
\truex{2500}%
\put(\value{x},0){\makebox(0,0)[l]{${#1}$}}%
\end{picture}}%

\newcommand{\Encurve}[2]%
{\begin{picture}(0,0)%
\truex{1400}\truey{2000}\truez{200}%
\put(\value{x},0){\oval(\value{y},#200)[r]}%
\put(\value{x},0){\makebox(0,0){\begin{picture}(0,#200)%
\put(0,0){\line(-1,0){\value{z}}}%
\put(0,#200){\vector(-1,0){\value{z}}}\end{picture}}}%
\truex{2500}%
\put(\value{x},0){\makebox(0,0)[l]{${#1}$}}%
\end{picture}}%

\newcommand{\Wscurve}[2]%
{\begin{picture}(0,0)%
\truex{1300}\truey{2000}\truez{200}%
\put(-\value{x},0){\oval(\value{y},#200)[l]}%
\put(-\value{x},0){\makebox(0,0){\begin{picture}(0,#200)%
\put(0,0){\vector(1,0){\value{z}}}%
\put(0,#200){\line(1,0){\value{z}}}\end{picture}}}%
\truex{2400}%
\put(-\value{x},0){\makebox(0,0)[r]{${#1}$}}%
\end{picture}}%

\newcommand{\Wncurve}[2]%
{\begin{picture}(0,0)%
\truex{1300}\truey{2000}\truez{200}%
\put(-\value{x},0){\oval(\value{y},#200)[l]}%
\put(-\value{x},0){\makebox(0,0){\begin{picture}(0,#200)%
\put(0,0){\line(1,0){\value{z}}}%
\put(0,#200){\vector(1,0){\value{z}}}\end{picture}}}%
\truex{2400}%
\put(-\value{x},0){\makebox(0,0)[r]{${#1}$}}%
\end{picture}}%

\newcount\SCALE%

\newcount\NUMBER%

\newcount\LINE%

\newcount\COLUMN%

\newcount\WIDTH%

\newcount\SOURCE%

\newcount\ARROW%

\newcount\TARGET%

\newcount\ARROWLENGTH%

\newcount\NUMBEROFDOTS%

\newcounter{x}%

\newcounter{y}%

\newcounter{z}%

\newcounter{horizontal}%

\newcounter{vertical}%

\newskip\itemlength%

\newskip\firstitem%

\newskip\seconditem%

\newcommand{\printarrow}{}%

\newcount\X

\newcommand{\truex}[1]{%
\NUMBER=#1%
\multiply\NUMBER by 100%
\divide\NUMBER by \SCALE%
\setcounter{x}{\NUMBER}}%

\newcommand{\truey}[1]{%
\NUMBER=#1%
\multiply\NUMBER by 100%
\divide\NUMBER by \SCALE%
\setcounter{y}{\NUMBER}}%

\newcommand{\truez}[1]{%
\NUMBER=#1%
\multiply\NUMBER by 100%
\divide\NUMBER by \SCALE%
\setcounter{z}{\NUMBER}}%

\newcommand{\changecounters}[1]{%
\SOURCE=\ARROW%
\ARROW=\TARGET%
\settowidth{\itemlength}{#1}%
\ifdim \itemlength > 2800\unitlength%
\addtolength{\itemlength}{-2800\unitlength}%
\TARGET=\itemlength%
\divide\TARGET by 1310%
\multiply\TARGET by 100%
\divide\TARGET by \SCALE%
\else%
\TARGET=0%
\fi%
\ARROWLENGTH=5000%
\advance\ARROWLENGTH by -\SOURCE%
\advance\ARROWLENGTH by -\TARGET%
\advance\SOURCE by -\TARGET}%

\newcommand{\initialize}[1]{%
\LINE=0%
\COLUMN=0%
\WIDTH=0%
\ARROW=0%
\TARGET=0%
\changecounters{#1}%
\renewcommand{\printarrow}{#1}%
\begin{center}%
\vspace{10pt}%
\begin{picture}(0,0)}%

\newcommand{\DIAGV}[2]{%
\SCALE=#1%
\setlength{\unitlength}{655sp}%
\multiply\unitlength by \SCALE%
\divide\unitlength by 100%
\initialize{\mbox{$#2$}}}%

\newcommand{\n}[1]{%
\changecounters{\mbox{$#1$}}%
\put(\COLUMN,\LINE){\makebox(0,0){\printarrow}}%
\thinlines%
\renewcommand{\printarrow}{\mbox{$#1$}}%
\advance\COLUMN by 4000}%

\newcommand{\nn}[1]{%
\put(\COLUMN,\LINE){\makebox(0,0){\printarrow}}%
\thinlines%
\ifnum \WIDTH < \COLUMN%
\WIDTH=\COLUMN%
\else%
\fi%
\advance\LINE by -4000%
\COLUMN=0%
\ARROW=0%
\TARGET=0%
\changecounters{\mbox{$#1$}}%
\renewcommand{\printarrow}{\mbox{$#1$}}}%

\newcommand{\conclude}{%
\put(\COLUMN,\LINE){\makebox(0,0){\printarrow}}%
\thinlines%
\ifnum \WIDTH < \COLUMN%
\WIDTH=\COLUMN%
\else%
\fi%
\setcounter{horizontal}{\WIDTH}%
\setcounter{vertical}{-\LINE}%
\end{picture}}%

\newcommand{\diag}{%
\conclude%
\raisebox{0pt}[0pt][\value{vertical}\unitlength]{}%
\hspace*{\value{horizontal}\unitlength}%
\vspace{10pt}%
\end{center}%
\setlength{\unitlength}{1pt}}%

\newcommand{\diagv}[3]{%
\conclude%
\NUMBER=#1%
\rule{0pt}{\NUMBER pt}%
\hspace*{-#2pt}%
\raisebox{0pt}[0pt][\value{vertical}\unitlength]{}%
\hspace*{\value{horizontal}\unitlength}%
\NUMBER=#3%
\advance\NUMBER by 10%
\vspace*{\NUMBER pt}%
\end{center}%
\setlength{\unitlength}{1pt}}%

\newcommand{\N}[1]%
{\raisebox{0pt}[7pt][0pt]{$#1$}}%

\newcommand{\crosslength}[2]{%
\settowidth{\firstitem}{#1}%
\settowidth{\seconditem}{#2}%
\ifdim\firstitem < \seconditem%
\itemlength=\seconditem%
\else%
\itemlength=\firstitem%
\fi%
\divide\itemlength by 2%
\hspace{\itemlength}}%

\usepackage{amsmath,amsthm,amsfonts,amssymb}

\usepackage{tocloft}

\usepackage{ifpdf}
\ifpdf       
\usepackage[pdftex, linktocpage]{hyperref}
\else
\usepackage[hypertex, linktocpage]{hyperref}
\fi
\usepackage{verbatim}

\title{O-minimal cohomology: finiteness and invariance results} 
\author{Alessandro Berarducci\thanks{Partially supported by the project:
Geometr\'{\i}a Real (GEOR) DGICYT MTM2005-02865 (2006-08)} \quad \& \ Antongiulio Fornasiero}
\date{May 26, 2007} 

\newlength{\tocwidth}
\makeatletter
\providecommand{\cftdotfill}{\@cftdotfill}
\makeatother

\makeatletter
\renewcommand{\@makecaption}[2]{%
\vskip\abovecaptionskip\sbox\@tempboxa{{\small#1. #2}}%
\ifdim\wd\@tempboxa>\hsize{{\small#1. #2}}\par%
\else\global\@minipagefalse\hb@xt@\hsize{%
  \hfil\box\@tempboxa\hfil}\fi\vskip\belowcaptionskip}
\makeatother

\DeclareMathOperator{\R}{\mathbb R}

\DeclareMathOperator{\Q}{\mathbb Q}
\DeclareMathOperator{\Z}{\mathbb Z}

\DeclareMathOperator{\im}{Im}

\DeclareMathOperator{\Ker}{Ker}

\DeclareMathOperator{\id}{id}
\DeclareMathOperator{\Sh}{Sh}
\newcommand{\ov}{\overline}

\newcommand{\nin}{\not\in}
\newcommand{\eps}{\varepsilon}
\newcommand{\ca}[1]{{\mathcal #1}}

\theoremstyle{plain}
\newtheorem{theorem}{Theorem}
\newtheorem{lemma}[theorem]{Lemma}
\newtheorem{proposition}[theorem]{Proposition}
\newtheorem{corollary}[theorem]{Corollary}
\newtheorem{conjecture}[theorem]{Conjecture}
\newtheorem{claim}{Claim}

\newtheorem{question}[theorem]{Question}
\newtheorem{fact}[theorem]{Fact}
\theoremstyle{definition}
\newtheorem{remark}[theorem]{Remark}
\newtheorem{definition}[theorem]{Definition}
\newtheorem{example}[theorem]{Example}
\newtheorem{exercise}[theorem]{Exercise}

\numberwithin{theorem}{section}
\numberwithin{equation}{section}

\newcommand{\bt}{\begin{theorem}}
\newcommand{\et}{\end{theorem}}

\newcommand{\bl}{\begin{lemma}}
\newcommand{\el}{\end{lemma}}

\newcommand{\bd}{\begin{definition}}
\newcommand{\ed}{\end{definition}}

\newcommand{\beq}{\begin{equation}}
\newcommand{\eeq}{\end{equation}}

\newcommand{\bexa}{\begin{example}}
\newcommand{\eexa}{\end{example}}

\newcommand{\bexe}{\begin{exercise}}
\newcommand{\eexe}{\end{exercise}}

\newcommand{\bfact}{\begin{fact}}
\newcommand{\efact}{\end{fact}}

\newcommand{\bprop}{\begin{proposition}}
\newcommand{\eprop}{\end{proposition}}

\newcommand{\bp}{\begin{proof}}
\newcommand{\ep}{\end{proof}}

\newcommand{\bc}{\begin{corollary}}
\newcommand{\ec}{\end{corollary}}

\newcommand{\bq}{\begin{question}}
\newcommand{\eq}{\end{question}}

\newcommand{\bcong}{\begin{conjecture}}
\newcommand{\econg}{\end{conjecture}}

\newcommand{\br}{\begin{remark}}
\newcommand{\er}{\end{remark}}

\newcommand{\Category}[1]{\mathfrak{#1}}
\newcommand{\Acat}{\Category{A}}
\newcommand{\Bcat}{\Category{B}}
\newcommand{\Ccat}{\Category{C}}
\DeclareMathOperator{\Xcat}{Sh_X}

\newcommand{\Fsh}{\ca F}
\newcommand*{\completion}[1]{\hat{#1}}
\newcommand{\Fshc}{\completion{\ca F}}
\newcommand{\Ab}{\Category{Ab}}
\newcommand{\Cov}[1]{\ca{#1}}
\newcommand{\Ucov}{\Cov{U}}
\newcommand{\Vcov}{\Cov{V}}

\DeclareMathOperator{\Der}{R}
\newcommand{\der}[1]{\Der^{#1}\!}
\newcommand{\dirlim}{\varinjlim}
\DeclareMathOperator{\Hom}{H}

\newcommand{\HomC}{\check{\Hom}\vphantom{\Hom}}
\newcommand{\rest}[1]{\mkern-3mu\upharpoonright\mkern-5mu%
\raisebox{-.5ex}{\ensuremath{#1}}}
\newcommand{\fc}{\check{f}}

\newcommand{\PCN}{PcN\xspace}

\newcommand{\cf}{cf.\xspace}
\newcommand{\ie}{i.e.\xspace}
\newcommand{\wloG}{w.l.o.g.\xspace}

\newcommand{\Cech}{\texorpdfstring{\v{C}ech}{Cech}\xspace}

\newcommand{\Pa}[1]{\bigl(#1\bigr)}
\newcommand{\set}[1]{\bigl\{#1 \bigr\}}

\newcommand{\abs}[1]{\lvert #1 \rvert}

\providecommand{\rom}[1]{{\/\rm #1}}

\newcommand{\w}{\widetilde}
\DeclareMathOperator{\basis}{{\cal U}}
\newcommand{\type}[1]{\langle #1 \rangle}

\begin{document}
\maketitle 

\vspace{-2em}
\begin{abstract} 
We prove that the cohomology groups of a definably compact set over an
o-minimal expansion of a group are finitely generated and invariant
under elementary extensions and expansions of the language. We also
study the cohomology of the intersection of a definable decreasing
family of definably compact sets under the
additional assumption that the o-minimal structure expands a field. 
\end{abstract} 

\begin{center}
\begin{minipage}{\tocwidth}
\small
\tableofcontents
\end{minipage}
\end{center}

\section{Introduction}
Delfs \cite{Delfs85} considered a sheaf cohomology theory for
(abstract) semialgebraic sets over arbitrary closed fields and proved
a semialgebraic version of homotopy invariance. In \cite{EdmundoJP05}
this was generalized to definable sets and maps over an o-minimal
expansion of a group. If one further assumes that the o-minimal
structure expands a field, then one can use the triangulation theorem
to show that the cohomology groups of a definable set are finitely
generated and invariant under both elementary extensions and
expansions of the language. We prove that this continues to hold for
arbitrary o-minimal expansions of groups, provided we restrict
ourselves to definably compact sets. 

Working without the field assumption entails various difficulties. 
To begin with one cannot make use of the apparatus of singular cohomology.
So we work, following the above authors, with sheaf cohomology. More
precisely, given a definable set $X\subset M^n$, the set of types $\w
X$ of $X$ with the ``spectral topology'' is quasi-compactification of
$X$, and we define $\Hom^i(X;\ca F ) := \Hom^i(\w X ; \ca F)$, where $\ca
F$ is a sheaf of Abelian groups on $\w X$.

Now consider the case when $G$ is an Abelian group and $\ca F$ is the
constant sheaf $G$ (\ie\ the sheaf generated by the presheaf with
constant value $G$). Assuming that $X$ is {\bf definably compact}
(\ie\ closed and bounded) we prove that $\Hom^i(X;G)$ is finitely
generated and invariant under both elementary extensions $N\succ M$
(\ie\ $\Hom^i(X;G) = \Hom^i(X(N);G)$) and expansions of the language of $M$
(note that expanding the language leaves $X$ invariant but alters $\w
X$).

This would be easy to prove if $M$ expands a field. In fact in this
case by the triangulation theorem $X$ is definably homeomorphic to the
geometrical realization $|K|$ (in $M$) of a finite simplicial complex
$K$, and a routine Mayer-Vietoris argument (together with the
acyclicity of simplexes) shows that $\Hom^i(X;G) \cong \Hom^i(K;G)$, where
the latter is the $i$-th simplicial cohomology group of $K$.

If $M$ does not expand a field we do not have the triangulation
theorem but we still have the cell decomposition theorem (see
\cite{Dries98}).  One could then be tempted to invoke the uniqueness
theorem for cohomology functors satisfying the (appropriate form) of
the Eilenberg-Steenrod axioms.  However, despite the cell
decomposition theorem, we have no guarantee that a definable set in an
o-minimal expansion of a group is a sort of definable CW-complex, so
the uniqueness theorem does not apply.  The problem is that in the
definition of a CW-complex one requires that the cells come equipped
with an attaching map that extends continuously to the boundary, while
for the o-minimal cells we do not have any such control of the
boundary.  

As standard references on sheaf cohomology and \Cech cohomology we
use \cite{Godement73} and \cite{Bredon97}. For the reader's
convenience we give in the appendixes the relevant definitions and
results.  Sheaf cohomology is defined for arbitrary topological
spaces, but many results are proved in the quoted texts under the
additional assumption that the space (or the family of supports) is
Hausdorff and paracompact. This is potentially a source of problems
since, given a definable set $X$, the spectral space $\w X$ associated
to it is in general not Hausdorff. In some cases we can reduce to the
compact Hausdorff situation using the fact that $\w X$ has a continuous
retraction onto a compact Hausdorff subspace ${\w X}^{\max}$ with the
same cohomology groups (see \cite{CarralC83}), but in some other cases
it is more convenient to show that the proofs of the relevant results
in \cite{Godement73} or \cite{Bredon97} work with the Hausdorff
hypothesis being replaced by normality (according to our convention
normality does not imply Hausdorff).  The latter approach has the
advantage that one can do without the spectrality hypothesis. In
particular it can be shown that for normal paracompact spaces (not
assumed to be Hausdorff or spectral) sheaf cohomology coincides with
\Cech cohomology with coefficients in the given sheaf. This is
reported in appendix~\ref{AP:D}, 
but we shall not really need it.
Instead we do need some results connecting the cohomology of
a subspace with the cohomology of its neighbourhoods.  Such results
have been proved in \cite{Delfs85,Jones06,EdmundoJP05} in the spectral
situation (see Corollary \ref{taut}), but they can also been
established under more general hypothesis (Theorem
\ref{LEM:TAUT-PCN}). As a corollary we obtain that 
$\Hom^*(\bigcap_{t>0} Y_t) \cong \dirlim_t \Hom^*(Y_t)$, whenever
$(Y_t \mid t>0)$ is a definable decreasing
family  of definably compact sets~$Y_t$
(Corollary~\ref{LEM:TAUT-DEFINABLE}).  We will also show
(Theorem~\ref{THM:DEFINABLE}) that if the o-minimal structure $M$
expands a field, then $\dirlim_t \Hom^*(Y_t) \cong \Hom^*(Y_{t_0})$ for
all sufficiently small $t_0$, but we are not able to prove this fact
without the field assumption.

\section{Topological preliminaries}\label{preliminaries} 

Let $X$ be a topological space. $X$ is {\bf normal} if every pair of
disjoint closed subsets of $X$ can be separated by open
neighbourhoods. $X$ is {\bf paracompact} if every open covering of $X$
has a locally finite refinement. Unlike other authors, in both
definitions we {\bf do not require that $X$ be Hausdorff}.  We shall
call a space {\bf \PCN} if it is paracompact and normal. Note that a
paracompact Hausdorff space is \PCN.

Let us recall that a {\bf quasi-compact} space is a topological space $X$ in
which every open covering has a finite refinement, or equivalently every
family of closed sets with the finite intersection property has a non-empty
intersection. So a compact space is an Hausdorff quasi-compact space.

Note that a quasi-compact space is \emph{a fortiori} paracompact.
Moreover, it is a well-known fact that a in a \PCN space $X$
every open covering $\Ucov$ admits a shrinking $\Vcov$:
\ie, $\Vcov$ is an open covering of~$X$,
and for every $V \in \Vcov$ there exists $U \in \Ucov$ containing~$\ov V$, the closure of~$V$.

A {\bf spectral space} is a
quasi-compact space having a basis of quasi-compact open sets stable
under finite intersections and such that every irreducible closed set
is the closure of a unique point. The prime spectrum of a commutative
ring with its Zariski topology is an example of a spectral
space. Another example is the set of prime filters of a lattice (see
\cite{CarralC83}). The set of $n$-types (ultrafilters of definable
sets) of a first order topological structure $M$ in the sense of
\cite{Pillay87} can also be endowed with a spectral topology (see
\cite{Pillay88}). In particular if $M$ is a real closed field one
obtains in this way the real spectrum of the polynomial ring $M[x_1,
\ldots, x_n]$ (see \cite{CosteR82}).

\section{Compactifications}

In this section we discuss a variant of the Wallman (or
Stone-\Cech) compactification of a normal topological space. The
variant depends on the particular choice of a basis of open sets. 

\bd \label{comp} 
Given a topological space $X$ with a fixed basis of open sets $\basis$
(that we always assume closed under finite intersections) we can
define a spectral space $\w X = \w {(X,\basis)}$ (depending on $\basis$) as follows.  A
{\bf constructible set} is a boolean combination of basic open sets
$U\in \basis$.  Let $\w X$ be the set of {\bf
ultrafilters} of constructible sets (\ie\ maximal families of
constructible sets closed under finite intersections and not
containing the empty set).  For $b\subset X$ constructible, let
$\tilde b = \{p\in \w X \mid b\in p\}$.  So $p\in \tilde b \Longleftrightarrow b\in
p$. The {\bf spectral topology} on $\w X$ is defined as follows. As a
basis of open sets of $\w X$ we take the sets of the form $\tilde b$
with $b$ an open constructible subset of $X$.
\ed

\bl $\w X$ is a spectral space. \el 

\bp To prove that $\w X$ is quasi-compact consider a family $\{C_i
\mid i\in I\}$ of closed sets $C_i \subset \w X$ with the finite
intersection property. We must prove that $\bigcap_{i\in I}C_i$ is
non-empty. Without loss of generality we can assume that $\{C_i\mid
i\in I\}$ is closed under finite intersections. Let $x\in \w X$ be an
ultrafilter containing all the closed constructible sets $b$ with $\w
b \supset C_i$ for some $i$. Then $x\in \bigcap_i C_i$, so $\w X$ is
quasi-compact. The same argument shows that the sets $\w b$ with $b$
constructible are quasi-compact.  So the sets $\w b$, with $b$ open
and constructible, form a basis of quasi-compact open sets of $\w X$
stable under finite intersections. To finish the proof we must show
that given an irreducible closed set $C$ of $\w X$, there is a unique
point $x\in C$ with $C = Cl(x)$.  To this aim, let $x$ be an
ultrafilter of constructible sets containing all the closed
constructible sets $b$ with $\w b \supset C$ and the complements of
the closed constructible sets $c$ such that $\w c \cap C$ is a proper
subset of $C$ (this family has the finite intersection property by the
irreducibility of $C$).  Then clearly $Cl(x) = C$.  To prove that $x$
is unique, suppose $Cl(x) = Cl(y)$. Then $x$ and $y$ contain the same
closed constructible sets. But the closed constructible sets generate
the boolean algebra of all the constructible sets. So $x$ and $y$ must
contain the same constructible sets, and are therefore equal. \ep

\br
Note that a point $x \in \w X$ is closed if and only if $x$
contains a maximal family of closed constructible sets with the finite intersection property. 
So by Zorn's lemma every closed subset $C$ of $\w X$ contains a closed point. 
\er

\noindent For $x\in X$ let \[\type x := \{b \mid x \in b\} \in \w X.\]
Since a constructible set $b$ is empty if and only if $\w b$ is empty,
the map \[\type {~} \colon X\to \w X\] has dense image.  Moreover this
map is injective whenever $X$ is a $T_1$-space. So in this case,
identifying $x$ with $\type x$, we have $X\subset \w X$, and it is
easy to see that the original topology on $X$ coincides with the
topology induced by $\w X$ (use the fact that for $A$ constructible,
$\w A \cap X = A$). Therefore:

\bl \label{T1} If $X$ is $T_1$, then for every open basis $\basis$ of~$X$,
$\w{(X, \basis)}$ is a quasi-compactification of~$X$. \el

We say that $(X,\basis)$ is {\bf constructibly
normal} if any pair of disjoint constructible open sets can be
separated by closed constructible sets. 

\bl If $(X,\basis)$ is constructibly normal, then $\w X$ is
normal (not necessarily Hausdorff).  \el

\bp Indeed given two disjoint closed subsets $A$ and $B$ of $\w X$, by
quasi-compactness there are disjoint closed constructible sets $A',B'$
in $X$ with $\w {A'} \supset A$ and $\w {B'} \supset B$. By the
assumption $A',B'$ can be separated by disjoint open constructible
sets $U\supset A'$ and $V \supset B'$. So $\w U$ and $\w V$ are open
sets separating $A,B$. 
\ep 

In a normal spectral space $Y$, the subset $Y^{\max}$ of the closed
points of $Y$ is compact Hausdorff (see \cite{CarralC83}). Also note
that, if $x$ is a closed point of $X$ and the singleton $\{x\}$ is
constructible, then $\type x$ is a closed point of $\w X$. So we have:

\bl If $(X,\basis)$ is constructibly normal and the points of $X$ are
closed and constructible, then ${\w X}^{\max}$ is a compactification
of $X$, namely it is a compact Hausdorff space containing $X$ as a
dense subspace.\el

If $X$ is a normal Hausdorff space, and we take as a basis of $X$ the family
$\basis$ of {\em all} its open subsets, then ${\w X}^{\max}$ is the
{\bf Wallman compactification} of $X$ \cite{Wallman38}. For
normal Hausdorff spaces it coincides with the {\bf Stone-\Cech
compactification} (see \cite[Thm. 3.6.22]{Engelking89}).

\bexa 
Consider the space $\Q$ with a basis $\basis$ of open sets given by the open intervals 
$(a,b)\subset \Q$, where we allow $a=-\infty$ or $b = +\infty$. Then 
\[\w {(\Q,\basis)} = \{a^-, a^+\}_{a\in \Q} \cup \R \cup \{\pm \infty\}\]
where: $a^+$ is the unique ultrafilter containing all sets of the form
$(a,b)$ with $b>a$; similarly $a^-$ contains all sets $(b,a)$ with
$b<a$; $+\infty$ is the ultrafilter containing $(a,\infty)$ for all $a\in \Q$; finally
$-\infty$ contains $(-\infty,a)$ with $a\in \Q$. We have $Cl(a^-) =
\{a^-,a\}$ and $Cl(a^+) = \{a^+,a\}$. The set of closed points is
\[\w {(\Q,\basis)}^{\max} = \R \cup \{\pm \infty\},\]
a two point-compactification of $\R$. 
\eexa

\bexa 
Consider the space $\Q$ with a basis $\basis$ of open sets given by the bounded open intervals 
$(a,b)\subset \Q$. Then 
\[\w {(\Q,\basis)}^{\max} = \R \cup \{\infty\},\]
the one-point compactification of $\R$. 
\eexa

\section{Cohomology of definable sets} 
Let $M = (M,<,\ldots)$ be an o-minimal structure (\cf~\cite{Dries98}) 
expanding a dense linear order $(M,<)$.
For instance $M$ may be a divisible ordered
group or a real closed field. We put on $M$ the topology generated by
the open intervals and on $M^n$ the product topology. If $X$ is a
subset of $M^n$ we put on $X$ the induced topology from $M^n$ unless
otherwise stated. By a {\bf definable} set we mean a first-order
definable (with parameters) subset of $M^n$ for some $n$. For instance
if $M$ is a real closed field, the definable sets are the
semialgebraic sets. Let $X\subset M^n$ be a definable set. Let
$\basis$ be the basis of $X$ consisting of the open definable subsets
of $X$. Define $\widetilde X = \w {(X,\basis)}$ as in Definition \ref{comp}.
Note that, by the o-minimality assumption, the definable sets
coincide with the constructible ones, namely every definable set is a
boolean combination of open definable sets. It follows that $\w X$ is
the set of types of $X$ (a type $p\in \widetilde X$ can be identified
with an ultrafilter of definable sets such that $X\in p$) endowed with
the spectral topology: as a basis of open sets we take the sets of
the form $\widetilde U$ with $U$ a definable open subset of $X$.
Since $X$ is Hausdorff, by Lemma \ref{T1} $\w X$ is a
quasi-compactification of $X$, in general not Hausdorff.

Let us now make the further assumption that $M$ is an o-minimal
expansion of an ordered group.  In this case one can use the $M$-valued
metric  $\abs{x-y}$ to show that
every definable set $X$
is {\bf definably normal} (any pair of disjoint definable closed sets
can be separated by definable open sets). Therefore in this case
$\widetilde X$ is {\bf normal} (not necessarily Hausdorff), and the
subspace ${\w X}^{\max}$ of its closed points is compact Hausdorff
(and contains $X$ as a topological subspace).

Given a definable set $X\subset M^n$ and a sheaf of Abelian groups
$\cal F$ on $\widetilde X$, we define, following
\cite{Delfs85,Jones06,EdmundoJP05}:
\begin{equation} \Hom^i(X;{\cal F} ):= \Hom^i(\widetilde X ; {\cal F}),
\end{equation}
see Definition~\ref{DEF:COM}.
Equivalently one can work with sheaves directly on $X$ rather than
$\widetilde X$ by considering $X$ not as a topological space but as a
site in the sense of Groethendieck (see \cite[\S 1.3]{CarralC83}). 
If $A$ is a definable subset of~$X$, we write $\Hom^*(A; \cal F)$  for
$\Hom^*(\widetilde A; \cal F \rest {\widetilde A})$.

If $f\colon X\to Y$ is a definable function, then $f$ induces a
function $\widetilde f \colon \widetilde X \to \widetilde Y$ by $f(p)
= \{Z \mid f^{-1}(Z) \in p\}$ where $Z$ ranges over the definable
subsets of $Y$ and we identify $p$ with an ultrafilter of definable
sets.  We have $\widetilde{f(X)} = \widetilde f (\widetilde X)$ and
${\widetilde f}^{-1}(\widetilde Z ) = \widetilde {f^{-1}(Z)}$. It
follows that if $f$ is continuous, $\widetilde f$ is continuous. So if
$f$ is an homeomorphism, then $\widetilde f$ is an homeomorphism.

If $G$ is an Abelian group and $X$ is a topological space, the
constant sheaf on $X$ with stalk $G$ will also be denoted $G$.  

We recall that a {\bf definable homotopy} between two definable functions
$f\colon X \to Y$ and $g\colon X\to Y$ is a definable continuous function
$F\colon I\times X\to Y$, where $I = [a,b]$ is some closed bounded interval in
$M$, such that $F(a,x) = f(x)$ and $F(b,x) = g(x)$ for every $x\in X$. 
Note that $F$ induces a map $\w F \colon \w {I \times X} \to \w Y$,
but in general $\w {I \times X} \neq \w I \times \w X$, so we cannot consider $\w F$
as a sort of ``homotopy'' parametrised by $\w I$. 
Nevertheless we have the following definable version of the homotopy axiom:

\bfact\label{homotopy}
{\rm(\cite{Delfs85,Jones06,EdmundoJP05})}  If $X,Y$ are definable
sets and $f,g\colon X \to Y$ are definably homotopic definable maps,
then $\widetilde f, \widetilde g \colon \widetilde X \to \widetilde Y$
(although not necessarily homotopic) induce the same homomorphism in
cohomology, namely for every Abelian group $G$ we have
\[f^* = g^* \colon \Hom^i(Y; G) \to \Hom^i(X;G).\] Similarly 
for maps $f,g\colon (X,A) \to (Y,B)$ of pairs.  \efact

This has been proved by \cite{Delfs85} in the semialgebraic case (\ie\ when $M$ is a real
closed field and $X,Y,f,g$ are semialgebraic). Jones \cite{Jones06}
extended it to the case of definable sets and maps in an o-minimal
expansion $M$ of an ordered field. In \cite{EdmundoJP05} it is shown
that it suffices that $M$ is an o-minimal expansion of an ordered
group. All proofs make use of the following lemma:

\bl {\rm(\cite[Prop.~4.7]{Delfs85})} Let $I\subset M$ be a closed and bounded interval $[a,b]$.
Then $\Hom^p(I; G) = 0$ for all $p>0$
and every Abelian group $G$. \el

Indeed it is easy to see the lemma holds for an arbitrary interval $I$, not necessarily closed and bounded
(note that in any case $\w I$ is quasi-compact). 

\section{Contractibility of cells}

Let $M$ be an o-minimal expansion of a group. 

\bl \label{contraction1} Let $I$ be a bounded interval in $M$ \rom(closed,
half-closed, or open\rom). Then $I$ is definably contractible to a
point. \el

\bp Let us consider the case $I = (a,b)$. Given $0< t \leq \frac{b-a}
2$, there is a (unique) definable continuous function $f_t \colon
(a,b) \to (a,b)$ such that $f_t$ is the identity on $[a+t,b-t]$ and it is
constant on both $[a,a+t]$ and $[b-t,b]$, with values $a+t$ and $b-t$
respectively (so the image of $f_t$ is $[a+t,b-t]$). Define $f_0 =
f$. Then $(f_t)_{0\leq t \leq \frac{b-a} 2}$ is a deformation retract
of $(a,b)$ to the point $\frac {a+b} 2$. The other cases are
similar. \ep

We will employ the following notation: given $B \subseteq M^{n-1}$ and \mbox{$f,g : B \to M$},
\[\begin{aligned}{}
(f,g)_B &:= \set{(x,y)\in M^{n-1} \times M: x \in B \ \& \ f(x) < y < g(x)},\\
[f,g]_B &:= \set{(x,y)\in M^{n-1} \times M: x \in B \ \& \ f(x) \leq y \leq g(x)},\\
\Gamma(f) &:= \set{(x,y)\in M^{n-1} \times M: x \in B \ \& \ f(x) = y},
\text{ the graph of } f.
\end{aligned}\]

\bl \label{contraction2} If $C$ is a bounded cell of dimension $m>0$
in $M^n$ then there is a deformation retract of $C$ onto a cell of
strictly lower dimension.  So by induction every bounded cell is
definably contractible to a point. \el

\bp If $C$ is the graph of a function we can reason by induction on
the dimension of the ambient space. So the only interesting case is
when $m>1$ and
$C = (f,g)_B$.
Let $h = \frac {f+g} 2$.
We will define
a deformation retract from $C$ to $\Gamma(h)$.  We can assume that
$h$ is a constant function, since we can reduce to this case by a
definable homeomorphism which fixes all but the last coordinate (just
take any constant function $h_1\colon B \to M$, and define $f_1,g_1$
so that they differ from $h_1$ by the same amount in which $f,g$
differ from $h$). Since $C$ is bounded, there are constants $a,b\in M$
such that $h$ is the constant function $\frac {a+b} 2$ and $(f,g)_B
\subset B\times (a,b)$.  By (the proof of) Lemma \ref{contraction1}
there is deformation retract of $(a,b)$ onto $\{\frac {a+b} 2\}$,
which induces a deformation retract of $B\times (a,b)$ onto $B \times
\{\frac {a+b} 2\}$, namely onto the graph of $h$. \ep

By Fact \ref{homotopy} we obtain:

\bc \label{acyclic-cell} If $C$ is a bounded cell of dimension $m$ in
$M^n$ then $\Hom^p(C;G) = 0$ for all $p>0$ and every Abelian group~$G$.
\ec

If we generalize slightly the definition of definable homotopy and allow the parameter of a homotopy to vary in the interval $[-\infty, + \infty]$, we get that Fact~\ref{homotopy} is still true, and therefore in Lemmata~\ref{contraction1} and~\ref{contraction2} we can drop the ``boundedness'' hypothesis.
Thus, Corollary~\ref{acyclic-cell} is true also for unbounded cells.

\section{Cells with non-acyclic closure} 

Let $M$ be an o-minimal expansion of a group and let $X\subset M^n$ be
a definable set. We will prove (Theorem \ref{fin-gen}) that the
cohomology groups $\Hom^p(X;G)$ of $X$ are finitely generated.  An
important special case is when $X$ is the closure $\ov C$ of a bounded
cell $C$. One may be tempted to conjecture that $\Hom^i(\ov C ; G)=0$ in
dimension $i>0$, but Example \ref{strange-cell} shows that in general
this is false.  Indeed similar examples show that $\Hom^1(\ov C ; G)$ can
have arbitrarily large finite rank.

\bexa \label{strange-cell} Let $M=(\R,<,+,\cdot)$. There is a bounded
cell $C$ of dimension $2$ in $\R^4$ whose closure $\ov C$ has an ``hole'',
namely it is definably homotopic to a circle. \eexa

\begin{floatingfigure}[r]{15ex}
\setlength{\unitlength}{3947sp}%
\begingroup\makeatletter\ifx\SetFigFont\undefined%
\gdef\SetFigFont#1#2#3#4#5{%
  \reset@font\fontsize{#1}{#2pt}%
  \fontfamily{#3}\fontseries{#4}\fontshape{#5}%
  \selectfont}%
\fi\endgroup%
\begin{picture}(770,766)(3443,-3144)
{\thinlines
\put(3826,-2761){\circle{750}}
}%
{\put(3826,-2761){\circle{336}}
}%
{\put(3996,-2741){\line( 1, 0){160}}
}%
\end{picture}
 
\caption{}
\label{annulus}
\end{floatingfigure}
Before giving the example, let us observe that
in Figure~\ref{annulus} (an annulus with a ray removed) we have
a space homeomorphic to a disk whose closure is homotopic to a circle.
However the space of Figure~\ref{annulus} is not a cell in
the sense of o-minimal cell decompositions. So we have to proceed
differently.  

\begin{figure}[htbp]
\begin{center}
\setlength{\unitlength}{3947sp}%
\begingroup\makeatletter\ifx\SetFigFont\undefined%
\gdef\SetFigFont#1#2#3#4#5{%
  \reset@font\fontsize{#1}{#2pt}%
  \fontfamily{#3}\fontseries{#4}\fontshape{#5}%
  \selectfont}%
\fi\endgroup%
\begin{picture}(5000,2700)(2300,-4148)
{\thinlines
\put(2401,-2011){\circle{1000}}
}%
{\put(4501,-2761){\line( 1, 0){1800}}
\put(6301,-2761){\line( 1,-2){600}}
\put(6901,-3961){\line(-1, 0){3000}}
\put(3901,-3961){\line( 1, 2){600}}
}%
{\put(4276,-3211){\line( 1, 0){2250}}
}%
{\put(4126,-3511){\line( 1, 0){2550}}
}%
{\put(4801,-2461){\line( 1, 0){1200}}
\put(6001,-2461){\line( 3, 2){900}}
\put(6901,-1861){\line(-2, 1){600}}
\put(6301,-1561){\line(-1, 0){1800}}
\put(4501,-1561){\line(-2,-1){600}}
\put(3901,-1861){\line( 3,-2){900}}
\put(4801,-2461){\line( 1, 1){600}}
\put(5401,-1861){\line( 1,-1){600}}
}%
{\put(3901,-1861){\line( 5, 1){1500}}
\put(5401,-1561){\line( 5,-1){1500}}
}%
{\multiput(5401,-1561)(0.00000,-120.00000){3}{\line( 0,-1){ 60.000}}
}%
{\put(5401,-2761){\line(-5,-4){1500}}
}%
{\put(5401,-2761){\line( 5,-4){1500}}
}%
{\put(4801,-3961){\line( 1, 2){600}}
}%
{\put(6001,-3961){\line(-1, 2){600}}
}%
{\put(4411,-2946){\line( 1, 0){1985}}
}%
{\put(5171,-1611){\line(-1, 0){775}}
}%
{\put(4836,-1681){\line(-1, 0){575}}
}%
{\put(4476,-1751){\line(-1, 0){365}}
}%
{\put(4471,-1751){\line( 6,-5){565.574}}
}%
{\put(4836,-1681){\line( 5,-6){340.984}}
}%
{\put(5176,-1616){\line( 2,-5){132.414}}
}%
{\put(5321,-1956){\line( 1, 0){170}}
}%
{\put(5181,-2091){\line( 1, 0){450}}
}%
{\put(5026,-2246){\line( 1, 0){755}}
}%
{\put(5626,-1611){\line( 1, 0){770}}
}%
{\put(5961,-1681){\line( 1, 0){580}}
}%
{\put(6346,-1756){\line( 1, 0){345}}
}%
{\put(5626,-1616){\line(-2,-5){135.172}}
}%
{\put(5966,-1686){\line(-1, 0){  5}}
\put(5961,-1686){\line(-5,-6){338.934}}
}%
{\put(6336,-1756){\line(-6,-5){562.623}}
}%
{\multiput(2401,-2011)(0.00000,90.00){4}{\line( 0, 1){ 60.00}}
}%
\put(3151,-1936){\makebox(0,0)[lb]{\smash{\SetFigFont{12}{14.4}{\rmdefault}{\mddefault}{\updefault}{$\approx$}%
}}}
\put(7126,-1936){\makebox(0,0)[lb]{\smash{\SetFigFont{12}{14.4}{\rmdefault}{\mddefault}{\updefault}{$=D$}%
}}}
\put(4051,-2611){\makebox(0,0)[lb]{\smash{\SetFigFont{12}{14.4}{\rmdefault}{\mddefault}{\updefault}{$\bigg\downarrow$}%
}}}
\end{picture} 
\caption{A two-dimensional cell in $\R^3$}
\label{cell}
\end{center}
\end{figure}
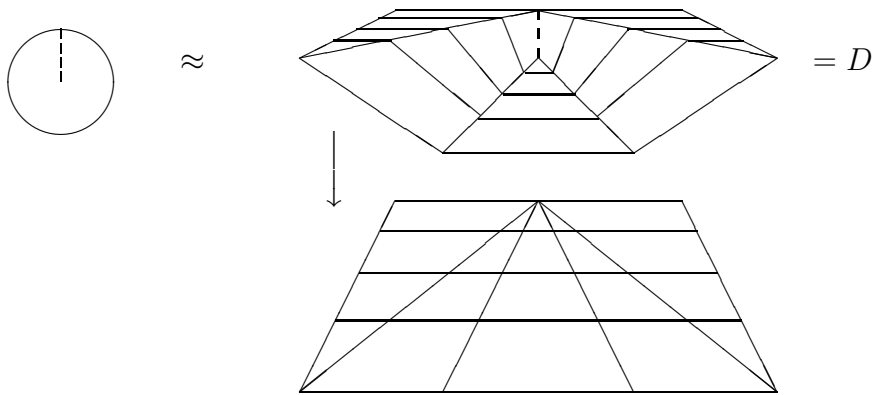

As a preliminary step we show that there is a
two-dimensional cell $D$ in $\R^3$ which is homeomorphic to an open
disk minus a ray via an homeomorphism which extends to the
closures. An example is the cell $D$ depicted on the top-right part of
Figure~\ref{cell}.

Figure~\ref{cell} must be interpreted as follows.  The projection of
$D$ on the first two coordinates is the open square on the
bottom-right part of Figure~\ref{cell}, with the five triangles
partitioning the square corresponding to the five regions in $D$. So $D$
is the graph of a function $f$ from the square to $\R$.  $D$ is
homeomorphic to an open disk with a ray removed via an homeomorphism
that extends to the closures and sends the ray of the disk to the
vertical dashed line $I \subset \partial D$ depicted in the Figure.
The three triangular regions of the cell $D$ are level sets of $f$,
with the central triangle being in a lower level, and the two top
triangles being on a common higher level.

We are now ready to describe the cell $C$ of Example
\ref{strange-cell}. The idea is the following.  Given a disk without a
ray as in Figure~\ref{cell}, we can ``bend'' it going in a higher
dimension, so as to create a ``hole'' (as in Figure~\ref{annulus}) by
separating the two sides of the middle part of the ray. More
precisely, consider a definable continuous function $g\colon D \to \R$
with the following properties: 1) $g$ takes non-zero values only for
points $x\in D$ sufficiently close to the middle point $x_0$ of the
dashed segment $I \subset \partial D$; 2) if $g(x) > 0$, then $x$ is
on the ``left hand side'' of $I$ and $g(x) \to 1$ for $x\to x_0$ from
the left; 3) if $g(x) < 0$, then $x$ is on the ``right-hand-side'' of
$I$ and $g(x) \to -1$ for $x\to x_0$ from the right. Let $C$ be the
graph of $g$. Then $H^1(\ov C; G) \neq 0$.

\br It is possible to show that $\Hom^1(X;G)$ can have arbitrarily high
rank even for a definable set $X$ that can be decomposed as a disjoint
union of only two bounded cells.  
Therefore, unlike the case of triangulations,
the abstract structure of a cell decomposition
(namely the dimensions and the adjacency relation between cells) by no
means determines the cohomology of a definable set. It is not however
excluded that by a refined version of the cell decomposition theorem one 
could obtain decompositions that do determine the cohomology groups.  \er

\section{Acyclic coverings}
In this section we give a sufficient condition, not based on
deformation retracts, to prove that an inclusion $X\subset Y$ of
topological spaces induces an isomorphism in cohomology (see Lemma
\ref{iso}).

Given an open cover $\ca U = \{U_i \mid i \in I\}$ of a topological
space $X$, and a subset $J\subset I$, we write $U_J$ for the
intersection $\bigcap_{i\in J} U_i$. The following theorem of Leray
says that, given an ``acyclic covering'', the cohomology of a sheaf
can be computed as the \Cech cohomology of the covering
(\cf~Appendix~\ref{AP:D}). 

\bfact\label{Leray}
{\rm(\cite[Thm. 4.13, p. 193]{Bredon97}, \cite[\S 5.4, p. 213]{Godement73})}
Let $\ca F$ be
a sheaf on a topological space $X$ and let $\ca U = \{U_i \mid i \in
I\}$ be an open covering of $X$ having the property that $\Hom^p(U_J; \ca
F) = 0$ for every $p>0$ and every finite $J\subset I$. 
Then the canonical homomorphism $\HomC^*(\ca U ; \ca F)\to \Hom^*(X;\ca F)$
is an isomorphism.
\efact

If moreover we assume that $\ca F$ is the constant sheaf $G$ and $\Hom^0(U_J; \ca F ) = G$ for every $J$ with $U_J$ non-empty, then from the
definition of \Cech cohomology it follows that $\HomC^*(\ca U;G)$
coincides with the $i$-th simplicial cohomology group with
coefficients in $G$ of the nerve $N({\ca U})$ of $\ca U$. So we have:

\bc \label{good-cover} 
If there exists a finite open cover $\ca U= \{U_i \mid i \in I\}$ of $X$ with 
 $\Hom^p(U_J; G ) = 0$ for every $p > 0$ and $\Hom^0(U_J; G) = G$ for every finite
$J\subset I$ with $U_J$ non-empty, then $\Hom^i(X;G)$ is isomorphic to the $i$-th 
simplicial cohomology group of the nerve $N({\ca U})$ of the covering, so in particular it is 
finitely generated. 
\ec 

Unfortunately we do not know the answer to the following:

\bq Let $X$ be a definably compact set in an o-minimal expansion $M$ of a
group. Does $\widetilde X$ have a cover as in Corollary \ref{good-cover}? \eq 

The answer is positive if $M$ expands a field, as a simple application
of the triangulation theorem shows.

\br \label{inducedh}
Let $f\colon X\to Y$ be a continuous function, let ${\cal V} =
\{V_i \mid i\in I\}$ be an (indexed) open cover of $Y$ and consider
the (indexed) open cover $f^{-1}({\cal V}):= \{ f^{-1}(V_i)\mid i \in I\}$ of
$X$. It follows easily from the definition of the induced
homomorphism $\fc$ in \Cech cohomology (see \cite[IX, \S 4]{EilenbergS52} or \cite[III.4.1.4]{Bredon97}) that there is a commutative
diagram:

\DIAGV{90} {\HomC^p({\cal V};G)} \n{\Ear{a}} \n{\HomC^p(f^{-1}({\cal V});G)} \nn 
           {\Sar{b}}        \n{}       \n{\Sar{c}} \nn 
           {\HomC^p(Y; G)}    \n{\Ear{\fc}} \n{\HomC^p(X; G)}        
\diag

\noindent where $b,c$ are the natural morphisms (see \cite{Godement73})
and $a$ is induced by the simplicial map on the nerves of
the indexed coverings sending $f^{-1}(V_i)$ to $V_i$.  \er

\bl \label{iso} Let $X\subset Y$ be definable sets. Suppose that there
are coverings ${\cal U} = \{ U_i \mid i\in I\}$ of $X$ and ${\cal V} =
\{ V_i \mid i \in I\}$ of $Y$ indexed by the same finite set $I$ such
that:
\begin{enumerate}
\item  $U_i \subset V_i$ for all $i,j\in I$.
\item For all finite $F\subset I$, $U_F$ is non-empty iff $V_F$ is
non-empty (\ie\ the natural map among the nerves of the coverings is
an isomorphism).
\item For each finite $F \subset I$ the sets $U_F$ and $V_F$ are either empty or connected, 
and for all $q>0$, $\Hom^q(U_F; G) = \Hom^q(V_F; G) = 0$.  
\end{enumerate}
Then the inclusion map $X\subset Y$ induces an isomorphism $\Hom^*(Y;
G) \to \Hom^*(X; G)$ for any Abelian group $G$. \el

Note that we do not require that $V_i \cap X = U_i$. 

\bp We are going to apply Remark \ref{inducedh} to the case when
$X\subset Y$ and $f$ is the inclusion map. In this case $f^{-1}({\cal
V}) = {\ca V} \cap X$ (by definition).  Consider the following
commutative diagram, where $a,b,c$ are as in Remark \ref{inducedh},
$d$ is induced by the natural simplicial isomorphism on the nerves of
the coverings, $e$ is the natural morphism form the
\Cech cohomology of a covering to the \Cech cohomology of the
space, and $p,q$ are the natural maps from \Cech to sheaf
cohomology.

\DIAGV{90} {\HomC^p({\cal V};G)} \n{\Ear{a}} \n{\HomC^p({\cal V}\cap X;G)} \n{\Ear{d}} \n{\HomC^p({\cal U}; G)}\nn           
           {\Sar{b}} \n{} \n{\Sar{c}} \n{\raisebox{1em}{\SWAR{\raisebox{1em}{e}}}} \n{} \n{} \nn
           {\HomC^p(Y; G)} \n{\Ear{\fc}} \n{\HomC^p(X; G)} \n{} \n{} \nn
           {\Sar{p}} \n{} \n{\Sar{q}} \n{} \n{} \nn
           {\Hom^p(Y; G)} \n{\Ear{f^* }} \n{\Hom^p(X; G)} \n{} \n{} 
\diag

\noindent By our assumptions on the coverings, $a$ and $d$ are isomorphisms. 
By Fact \ref{Leray} $p\circ b$ and $q\circ e$ are isomorphisms. So 
 $f^*: \Hom^*(Y;G) \to \Hom^*(X; G)$  is
an isomorphism.  \ep

\section{Finiteness results for cohomology} 

As usual let $M$ be an o-minimal expansion of group. Given a definably
compact set $X$ we want to prove that $\Hom^i(X;G)$ is finitely generated
for every $i$ and every Abelian group $G$. An important special case
is when $X$ is the closure $\ov C$ of a bounded cell $C$. We will show
(as a consequence of Corollary \ref{X-C}) that there is a point $a\in
C$ such that $\Hom^p(\ov C \setminus \{a\};G) \cong \Hom^p(\partial C;G)$,
where $\partial C := \ov C \setminus C$ is the boundary of
$C$. Granted this, since $\partial C$ has smaller dimension than $\ov
C$, by induction on the dimension we can assume that the cohomology
groups of $\partial C$ are finitely generated and carry on with the
inductive proof. At first sight the fact that $\Hom^p(\ov C \setminus
\{a\};G) \cong \Hom^p(\partial C;G)$ looks rather intuitive: one may even
be tempted to conjecture that $\partial C$ is a definable deformation
retract of $\ov C \setminus \{a\}$ (as it would be the case for the
cells of a CW-complex), or at least that these two spaces are
definably homotopy equivalent. However we are not able to prove this.
We proceed instead in a different manner, with the role of definable
deformation retracts being taken by Lemma \ref{iso}.  There is however
a further complication. We are not able to apply Lemma~\ref{iso}
directly to the pair of sets $(X,Y) = (\ov C \setminus \{a\}, \partial
C)$, but only to pairs of sets of the form $(C \setminus \{a\},
C\setminus C_t)$, where $(C_t)_{t>0}$ is a suitable definable
collection of sets with $\bigcup_{t>0} C_t = C$, and consequently $\bigcap_{t>0}(\ov C \setminus C_t) = \partial C$ (the singleton $\{a\}$ is one of the
$C_t$). So we first prove $\Hom^p(C \setminus \{a\}) \cong \Hom^p(C\setminus
C_t)$. Then we deduce $\Hom^p(\ov C \setminus \{a\}) \cong \Hom^p(\ov
C\setminus C_t)$ by the excision theorem. Finally, with the help of
Corollary \ref{LEM:TAUT-DEFINABLE}, we let $t\to 0$ to obtain $\Hom^p(\ov
C \setminus \{a\}) \cong \Hom^p(\partial C)$. This is the idea. Let us
now come to the details.

\bl \label{cell-covers}
Let $C\subset M^n$ be a bounded cell of dimension $m$.
There is a definable family $\{ C_t \mid t>0 \}$ of definably compact sets
$C_t\subset C$ such that:
\begin{enumerate}
\item $C = \bigcup_{t>0}C_t$.
\item If $0<t'<t$, then $C_t\subset C_{t'}$ and the inclusion
$C\setminus C_{t'} \subset C\setminus C_{t}$ induces an isomorphism
$\Hom^p(C \setminus C_t; G) \to \Hom^p(C \setminus C_{t'}; G)$.
\item For every Abelian group~$G$,
$C\setminus C_t$ has the same cohomology groups of an $m-1$
dimensional sphere, namely $\Hom^p(C\setminus C_t; G) = 0$ for $p \nin
\{0,m-1\}$, and $\Hom^p(C\setminus C_t; G) = G$ for $p \in \{0,m-1\}$.
\end{enumerate}
\el

Note that the result would be obvious if $M$ expands a field, since in
that case $C$ is definably homeomorphic to an open ball, and each open
ball is the increasing union of its concentric closed sub-balls.

\noindent\textit{Proof}.
We define $C_t$ as follows. 

\begin{enumerate}
\item If $n=1$ and $C = (a,b)$, then $C_t = [a+ \gamma_t, b-
\gamma_t]$ where $\gamma_t = \min{\{\frac {a+b} 2, t \}}$ (so $C_t$ is
non-empty).
\item If $n=1$ and $C$ is a singleton in $M$, $C_t = C$.
\item Let $n>1$ and $C = \Gamma(f)$, where $f\colon B\to M$. By
induction $B_t$ is defined and we set $C_t = \Gamma(f \rest{B_t})$. 
\item 
Let $n>1$ and $C = (f,g)_B$. By induction $B_t$ is defined.
We put
$C_t = [f+\gamma_t, g-\gamma_t]_{B_t}$,
where $\gamma_t := \min(\frac{f - g}{2}, t)$.
\end{enumerate} 

With this definition we have:

\begin{claim} \label{cell-cover} For each $t>0$ 
there is a covering $\ca U = \{U_i\mid i\in I\}$ of
$C\setminus C_t$ such that:
\begin{enumerate}
\item 
The index set $I$ is the family of the closed faces of an
$m$-dimensional cube, where $m= \dim(C)$. \rom(So $|I|=2m$\rom).
\item If $F\subset I$, then $U_F := \bigcap_{i\in F}U_i$ is either
empty or a cell. \rom(So in particular $\Hom^p(U_F; G) = 0$ for all $p>0$
and, if $U_F \neq \emptyset$, $\Hom^0(U_F;G) = G$.\rom)
\item For $F\subset I$, $U_F \neq \emptyset$ iff the faces of the
cubes belonging to $F$ have a non-empty intersections. \rom(So the nerve of
$\ca U$ is isomorphic to the nerve of a covering of an $m$-cube by its
closed faces.\rom)
\end{enumerate} 
\end{claim} 
\begin{minipage}{\textwidth}
\begin{floatingfigure}[r]{25ex}
\setlength{\unitlength}{2565sp}%
\begingroup\makeatletter\ifx\SetFigFont\undefined%
\gdef\SetFigFont#1#2#3#4#5{%
  \reset@font\fontsize{#1}{#2pt}%
  \fontfamily{#3}\fontseries{#4}\fontshape{#5}%
  \selectfont}%
\fi\endgroup%
\begin{picture}(2424,2124)(3889,-4273)
\thinlines
{\put(4006,-3856){\oval(210,210)[bl]}
\put(4006,-2566){\oval(210,210)[tl]}
\put(4396,-3856){\oval(210,210)[br]}
\put(4396,-2566){\oval(210,210)[tr]}
\put(4006,-3961){\line( 1, 0){390}}
\put(4006,-2461){\line( 1, 0){390}}
\put(3901,-3856){\line( 0, 1){1290}}
\put(4501,-3856){\line( 0, 1){1290}}
}%
{\put(4306,-2656){\oval(210,210)[bl]}
\put(4306,-2266){\oval(210,210)[tl]}
\put(5896,-2656){\oval(210,210)[br]}
\put(5896,-2266){\oval(210,210)[tr]}
\put(4306,-2761){\line( 1, 0){1590}}
\put(4306,-2161){\line( 1, 0){1590}}
\put(4201,-2656){\line( 0, 1){390}}
\put(6001,-2656){\line( 0, 1){390}}
}%
{\put(5806,-3856){\oval(210,210)[bl]}
\put(5806,-2566){\oval(210,210)[tl]}
\put(6196,-3856){\oval(210,210)[br]}
\put(6196,-2566){\oval(210,210)[tr]}
\put(5806,-3961){\line( 1, 0){390}}
\put(5806,-2461){\line( 1, 0){390}}
\put(5701,-3856){\line( 0, 1){1290}}
\put(6301,-3856){\line( 0, 1){1290}}
}%
{\put(4306,-4156){\oval(210,210)[bl]}
\put(4306,-3766){\oval(210,210)[tl]}
\put(5896,-4156){\oval(210,210)[br]}
\put(5896,-3766){\oval(210,210)[tr]}
\put(4306,-4261){\line( 1, 0){1590}}
\put(4306,-3661){\line( 1, 0){1590}}
\put(4201,-4156){\line( 0, 1){390}}
\put(6001,-4156){\line( 0, 1){390}}
}%
\end{picture}
 
\caption{}
\label{cover}
\end{floatingfigure}
For example for $m=2$ we have four open sets $U_i$ which intersect each other
as in Figure~\ref{cover}. 
Note that the claim implies, by Corollary~\ref{good-cover}, that $C
\setminus C_t$ has the same cohomology groups of an $m-1$ dimensional
sphere.  To prove the claim we define $\ca U$ by induction on the
dimension $n$ of the ambient space. We distinguish four cases
according to the definition of~$C_t$.
\end{minipage}
\begin{enumerate}
\item{If $n=1$ and $C$ is a singleton, then $\ca U$ is the covering
consisting of one open set (given by the whole space $C$).}
\item{If $n=1$ and $C = (a,b)$, then $C\setminus C_t$ is the union of
the two open subsets $(a, a + \gamma_t)$ and $(b- \gamma_t, b)$,
and we define $\ca U$ as the covering consisting of these two sets.}
\item Let $n>1$ and $C = \Gamma(f)$, where $f\colon B\to M$. By
induction we have a covering $\ca V$ of $B\setminus B_t$ with the
stated properties, and we define $\ca U$ to be the covering of
$C\setminus C_t$ induced by the natural homeomorphism between the
graph of $f$ and its domain.
\item
Let $n > 1$ and $C = (f,g)_B$. By definition $C_t = (f + \gamma_t, g -
\gamma_t)_{B_t}$.  By induction $B\setminus B_t$ has a covering $\ca V
= \{V_j \mid j\in J\}$ with the stated properties, where $J$ is the
set of closed faces of the cube $[0,1]^{m-1}$.  Define a covering $\ca
U = \{U_i \mid i\in I\}$ of $C\setminus C_t$ as follows. As index set
$I$ we take the closed faces of the cube $[0,1]^m$. Thus $|I| = |J| +
2$, the two extra faces corresponding to the ``top'' and ``bottom''
face of $[0,1]^m$. 
We associate to the top face 
the open set $(g- \gamma_t,g)_B \subset C\setminus C_t$ and to the
bottom face the open set $(f, f + \gamma_t)_B \subset C \setminus
C_t$. The other open sets of the covering are the preimages of the sets $V_j$ 
under the projection $M^n\to M^{n-1}$.  This
defines a covering of $C\setminus C_t$ with the stated properties.
\end{enumerate}

It remains to show that the inclusion map $C\setminus C_{t'} \subset
C\setminus C_{t}$ induces an isomorphism $\Hom^p(C \setminus C_t; G) \to
\Hom^p(C \setminus C_{t'}; G)$.  To this aim it suffices to observe that
by (the proof of) Claim 1 there are coverings $\ca U$ of $C \setminus
C_{t'}$ and $\ca V$ of $C \setminus C_t$ satisfying the assumptions of
Lemma \ref{iso}. 
\hfill\textsquare

\bl \label{x-c} Let $X$ be a definably compact set, $C$~be a cell of
maximal dimension in~$X$, and $G$ be an Abelian group.
Then for each $0<t'<t$ the inclusion map
$X\setminus C_{t'} \subset X\setminus C_t$
induces an isomorphism
\[\Hom^*(X\setminus C_{t}; G) \cong \Hom^*(X\setminus C_{t'}; G).\] \el

\bp 
By the excision theorem  (\cite[Thm. 12.9]{Bredon97})
the inclusion of pairs $(C\setminus C_{t'}, C\setminus C_t) \to (X
\setminus C_{t'}, X \setminus C_{t})$ induces an isomorphism \[\Hom^*(X
\setminus C_{t'}, X \setminus C_{t}; G) \cong \Hom^*(C\setminus C_{t'},
C\setminus C_t; G).\] The right-hand side is zero by
Lemma \ref{cell-covers} and the long cohomology sequence of the pair (see
\cite[eqn. 24, p. 88]{Bredon97}). So the left-hand side is also zero and therefore
$\Hom^*(X\setminus C_{t}; G) \to \Hom^*(X\setminus C_{t'}; G)$ is an isomorphism. \ep

In the above theorem it is not necessary that $X$ is definably
compact: it suffices that $C$ is bounded.

\br \label{lim} Let $R$ be the limit of a direct system
$(R_i \mid f_{i,j})_{i,j\in I}$ of Abelian groups and morphisms
$f_{i,j}\colon R_i \to R_j$. 
Suppose that each $f_{i,j}$ is an isomorphism. Then for each
$i$ the natural morphism $R_i \to R$ is an isomorphism. \er

\bc \label{X-C} Suppose that $X$ is a definably compact set such that
$C$ is a cell of maximal dimension of $X$ \rom(for instance $X = Cl(C)$\rom).
Then for every $t>0$ the inclusion map
$X\setminus C \subset X\setminus C_t$ induces an isomorphism \[\Hom^p(X
\setminus C_t; G) \cong \Hom^p(X \setminus C; G).\]  \ec

\bp  For $t>0$ we have
\[ \begin{array}{ccc}
\Hom^*(X\setminus C_t;G) & \cong & \dirlim_s \Hom^*(X \setminus C_s;G) \\
 & \cong & \Hom^*(X \setminus C;G)
\end{array} \]
where the first isomorphism follows from Lemma~\ref{x-c} and
Remark~\ref{lim}, and the second one follows from the equality $X
\setminus C = \bigcap_s (X \setminus C_s)$ and
Corollary~\ref{LEM:TAUT-DEFINABLE}. \ep

Note that for $t$ big enough, $C_t$ is a singleton. So in particular,
applying the theorem to $X = \ov{C}$, we have proved that there is a
point $a\in C$ such that:
 \begin{equation} \label{partial} \Hom^p(\ov{C}
\setminus \{a\}; G) \cong \Hom^p(\partial C; G). \end{equation} We do not
know however whether $\partial C$ is a definable deformation retract of $\ov C
\setminus \{a\}$.

\bt \label{fin-gen} Let $X \subset M^n$ be a definably compact
set and $G$ be an Abelian group.
Then, for each $p$, $\Hom^p(X;G)$ is finitely generated. Moreover
$\Hom^p(X;G)=0$ for $p>\dim(X)$.  \et

\begin{proof}
By a Mayer-Vietoris argument. Decompose $X$ into cells. Let $C$ be
a cell of $X$ of maximal dimension. Let $t>0$ be such that $C_t \neq
\emptyset$ and write $X$ as the union of $X\setminus C_t$ and
$C$. Consider the Mayer-Vietoris sequence (see \cite[eqn. 32,
p. 98]{Bredon97}) associated to this union:
\begin{equation} \label{EQ:MAYER}
\ldots \to \Hom^{p-1}(C \setminus C_t) \to \Hom^p(X) \to \Hom^p(X\setminus C_t)
\oplus \Hom^p(C) \to \Hom^p(C\setminus C_t) \to \ldots
\end{equation}
where we have omitted the coefficients $G$ in the notation. By Corollary \ref{X-C} the
inclusion $X\setminus C_t \subset X\setminus C$ induces an isomorphism
in cohomology, so composing with this isomorphism we obtain 
\begin{equation} \label{mv} 
\ldots \to \Hom^{p-1}(C \setminus C_t) \to \Hom^p(X) \to \Hom^p(X\setminus C) \oplus
\Hom^p(C) \to \Hom^p(C\setminus C_t) \to \ldots
\end{equation}
Now $C$ has the cohomology of a point and $C\setminus C_t$ has the
cohomology of an $(m-1)$-dimensional sphere. So the displayed part of
the sequence above has the form:
\begin{align}\label{mv1} 
& G \to \Hom^p(X) \to \Hom^p(X\setminus C) \oplus
\Hom^p(C) \to 0,\\
\text{ or}\quad &
\label{mv2} 
0 \to \Hom^p(X) \to \Hom^p(X\setminus C) \oplus
\Hom^p(C) \to G,\\
\text{ or} \quad &
\label{mv3} 
0 \to \Hom^p(X) \to \Hom^p(X\setminus C) \oplus
\Hom^p(C) \to 0;
\end{align}
where (\ref{mv3}) applies for $p \in \{m, m-1\}$. 
 By induction on the number of cells 
$\Hom^p(X\setminus C)$ is finitely generated, and vanishes for $p\geq m$.
 From the above sequences it then follows that the same holds for $\Hom^p(X)$.
\end{proof}

\section{Elementary extensions and change of language}

Let $M$ be an o-minimal expansion of a group and let $X\subset M^n$ be
a definable set.  We have seen that we can associate to $X$ the
spectral space $\w X$ of all its types over $M$ (such a type can
be identified with an ultrafilter of $M$-definable sets that contains 
$X$). The cohomology of $X$ has been defined as the cohomology of $\w X$. 
If $N\succ M$ is an elementary extension we may also associate to $X$ the spectral
space $\w{X(N)}$ (ultrafilters of $N$-definable sets that contains $X(N)$). 

\bt \label{theta} Let $\theta \colon \widetilde{X(N)} \to \w X$ be
the map that sends a type over $N$ to its restriction over $M$. If $X$
is definably compact then $\theta$ induces an isomorphism
$\Hom^*(\w X; G) \to \Hom^*(\widetilde{X(N)}; G)$ for any Abelian group $G$. \et

\bp
For both $X$ and $X(N)$ we have an exact sequence as in
equation~\eqref{EQ:MAYER}
above. The terms of the two exact sequences so obtained are connected by the
homomorphisms induced by $\theta$ in such a way that the resulting 
diagram commutes (it is important that in equation~\eqref{EQ:MAYER}
we take the parameter $t$ in the small model~$M$).
The desired result then follows
arguing as in Theorem \ref{fin-gen} by induction on the number of cells. \ep

We now extend the above result to certain type-definable sets. As above let $N\succ M$. 

\bt Let $X\subset M^n$ be a definably compact set. Let $A\subset
\w X$ be a type-definable closed subset and let $A(N) := \theta^{-1}(A)
\subset \widetilde{X(N)}$. Then $\theta$ induces an isomorphism
$\Hom^*(A; G) \to \Hom^*(A(N) ; G)$. In particular if
$p$ is a closed type in $\w X$, then the set $\theta^{-1}(p)$ of all
types of $\widetilde{X(N)}$ which restrict to $p$ has the same
cohomology of a point \rom(so in particular it is connected\rom). \et
\bp Each closed type-definable set $A\subset \w X$ can be written as an
intersection $\bigcap_{i\in I}X_i$ of definably compact sets $X_i
\subset M^n$. Now observe that $A(N) =
\bigcap_{i\in I}\widetilde{X_i(N)}$.  By Corollary \ref{LEM:TAUT-DEFINABLE} and Theorem
\ref{theta} $\Hom^*(A(N);G) = \dirlim_{i\in I}\Hom^*(\widetilde{X_i(N)};G) =
\dirlim_{i\in I}\Hom^*(X_i;G) = \Hom^*(A;G)$.
\ep
In similar way we can prove: 

\bt If $M_1$ is an o-minimal expansion of $M$ to a bigger language and
$X\subset {M}^n$ is a definably compact set in $M$, then the map
$\widetilde{X(M_1)} \to \w X$ sending each type in the language
$L_1$ to its restriction to $L$ induces an isomorphism $\Hom^*(\w
X; G) \to \Hom^*(\widetilde{X(M_1)}; G)$ for any Abelian group $G$.\et

\section{Definable families in expansions of fields}
\begin{definition}
Given definable maps $f : X \to Y$ and $f' : X' \to Y'$, we say that $f$ and $f'$ have the same definable topological type iff there exist definable homeomorphisms $\lambda: X \to X'$ and $\mu: Y \to Y'$ making the following diagram commute:
\DIAGV{50}
{}\n{X} \n{\Ear{f}} \n{X'} \nn
{}\n{\Sar{\lambda}} \n{} \n{\saR{\mu}} \nn
{}\n{Y} \n{\Ear{g}} \n{Y'}
\diag
\end{definition}

\bd
Let $(Y_t)_{t > 0}$ be a family of subsets of~$M^n$.
We will say that $(Y_t)_{t>0}$ is decreasing if $Y_t \subseteq Y_{t'}$ for every $t \leq t'$.
\ed

If $M$ is an o-minimal expansion of a
field, Corollary \ref{LEM:TAUT-DEFINABLE} can be strengthened as follows. 

\bt \label{THM:DEFINABLE}
Assume that $M$ is an o-minimal expansion of a field.
Let $\Pa{Y_t}_{t > 0}$ be a definable decreasing
family of definably compact subsets of some definable set~$Y$.
Let $A := \bigcap_{t > 0} Y_t$. 
Then for every sufficiently small $t$ we have a natural isomorphism
induced by the inclusion
\[
\Hom^*\Pa{ A; {\ca F}} \cong \Hom^*\Pa{{Y_t}; {\ca F}},
\]
for every sheaf ${\ca F}$ on~$\widetilde Y$.
\et
In the above theorem we cannot weaken the hypothesis to $Y_t$ closed
(instead of definably compact).  For instance, let $Y_t :=
[\frac{1}{t}, \infty) \subseteq M$. 

Note that Corollary \ref{X-C} is a special case of Theorem \ref{THM:DEFINABLE}.
To prove the theorem, we need the following lemmata.
Note that for the lemmata we do not need that $M$ expands a field,
but only that it expands a group.
\bl \label{REM:MINIMA} Let $M$ be an o-minimal expansion of a group. 
Let $Y \subseteq M^n$ be definable, and $f: Y \to M$ be a definable function \rom(not necessarily continuous\rom).
Let $A \subseteq Y$ be the set of local minima of~$f$.
Then, $f(A)$ is finite.
\el 
\begin{proof}
If not, let $a < b \in M$ such that, for every $t \in (a,b)$ there
exists $\gamma(t) \in Y$ such that $f(\gamma(t)) = t$ and $\gamma(t)$
is a local minimum for~$f$.  By definable choice, we can assume that
$\gamma$ is a definable continuous function. It follows that in
 any neighbourhood of $\gamma(t)$ there are points of the
form $\gamma(t')$ with $t'<t$. But $f(\gamma(t')) < f(\gamma(t))$, contradicting the fact that $\gamma(t)$ is a local minimum.
\end{proof}

\bl\label{LEM:FAMILY}
Let $M$ be an o-minimal expansion of a group. 
Let $(Y_t)_{t > 0}$ be a decreasing definable family of closed subsets of~$M^n$.
Then, there exists $t_0 >0$ such that, for every $t \in M$ with $0 < t < t_0$ we have
\[
Y_t = \bigcap_{t < u < t_0} Y_u
\]
\el
\begin{proof}
Let $Y := \bigcup_{t > 0} Y_t$.
Define
\begin{align*}
f &: Y \to M^{\geq 0}\\
f(x) &:= \inf \set{t: x \in Y_t}.
\end{align*}
Define also $Z_t := f^{-1}([0,t])$.
Note that $Z_t = \bigcap_{t < u < t_0} Y_u$, and
$Y_t \subseteq Z_t$.
Therefore,
the conclusion is equivalent to saying that, for every $0 < t < t_0$, we have that $Y_t = Z_t$.

\begin{claim}
For every $t \geq 0$, 
if $x \in Z_t \setminus Y_t$, then $x$ is a local minimum for~$f$.
\end{claim}
If not, let  $x \in Z_t \setminus Y_t$ such that $x$ is not a local minimum for~$f$. Note that $f(x) = t$.
Let $\gamma: (0, \eps) \to Y$ be a definable function such that $\lim_{s \to 0}\gamma_s = x$ and $f(\gamma_s) < f(x)$.
Then, $\gamma_s \in Y_t$.
Since $Y_t$ is closed, we have that $x \in Y_t$, absurd.

By Lemma~\ref{REM:MINIMA}, there exists $t_0 >0$ such that, for every 
$x \in Y$, if $0 < f(x) < t_0$, then $x$ is not a local minimum for~$f$.
The claim implies the conclusion.
\end{proof}

\br
A few remarks about the above Lemma and its proof.
\begin{enumerate}
\item The function $f$ is lower semi-continuous, because $Z_t$ is closed for every~$t$.
\item The hypothesis that the $Y_t$ are closed is necessary.
\item Even if the function $f$ in the proof is continuous, we cannot conclude that $Y_t = Z_t$ always.
\item It is not true that there exists $t_0 > 0$ such that $f$ is continuous on~$Y_{t_0}$.
\end{enumerate}
\er 

\begin{proof}[Proof of Theorem~\ref{THM:DEFINABLE}]
By Corollary~\ref{LEM:TAUT-DEFINABLE},
\[
\Hom^*(\w A; {\ca F}) = \dirlim_{t \to 0} \Hom^*(\w{Y_t}; {\ca F}).
\]
Since $\Pa{Y_t}_{t > 0}$ is a definable family, and $M$ expands a
field, the (definable) topological type of the $Y_t$ is eventually
constant (this is a consequence of the trivialization theorem holding
in o-minimal expansions of fields, see \cite{Dries98}). Hence we can
assume \wloG\ that the $Y_t$ are all definably homeomorphic, and
therefore that their cohomology groups over $G$ are isomorphic to each
other; let us denote by $P$ such group.  Thus, for every $t \leq t'$,
the map induced by the inclusion
\[
\phi^{t'}_{t}: \Hom^*(\w{Y}_{t'}; {\ca F}) \to \Hom^*(\w{Y}_{t}; {\ca F})
\]
is an endomorphism of~$P$.
It suffices to prove that $\phi^{t'}_{t}$ is an isomorphism for all sufficiently small $t'$ to prove the theorem.
By Lemma~\ref{LEM:FAMILY}, \wloG we can assume that, for every $t >0$,
\[
Y_t = \bigcap_{t' >t} Y_{t'}.
\]

For convenience, define $\phi^{t'}_t := \phi^{t}_{t'}$  if $t' \leq t$, and $I := M^{>0}$.
Let $E$ be the equivalence relation on $I$ given by
\[
t E t' \Leftrightarrow \phi^{t'}_{t} \text{is an isomorphism}.
\]
\begin{claim}
$E$~is a definable subset of $I^2$.
\end{claim}
In fact, by the trivialization theorem, there exists a finite
partition $\{ C_1, \dotsc, C_n \}$ of $\{(t,t') \in M^2: 0 < t \leq
t' \}$ into definable sets, such that the definable topological type of
the inclusion map $Y_t \to Y_{t'}$ is constant on each~$C_k$.
Hence, the isomorphism type of the map $\phi^{t'}_t$ is constant on each
$C_k$, and thus $E$ is a finite union of some of the~$C_k$.

\begin{claim}\label{CL:OPEN}
For every $t \in I$ there exists $t_0 >t$ such that, for every $t' \in I$,
\[
t \leq t' < t_0 \Rightarrow t E t'.
\]
\end{claim}\noindent
Since $Y_t = \bigcap_{t' >t} Y_{t'}$, by Corollary
\ref{LEM:TAUT-DEFINABLE} we have $\Hom^*(\w{Y_t}; {\ca F}) =$%
\linebreak[2]%
\raisebox{0pt}[2ex][0pt]{$\displaystyle{\dirlim_{t' \to t^+}} \Hom^*(\w{Y}_{t'}; {\ca F})$};
namely,
\[
P = \dirlim_{t' \to t^+} \set{P, \phi^{t'}_{t}}.
\]
By definition of inductive limit, for every $p \in P$ there exists $t'>t$ and $p' \in P$ such that $p = \phi^{t'}_{t}(p')$.
However, $P$ is finitely generated, and therefore there exists $t_0 >t$ such that $\phi^{t'}_t$ is surjective for every $t_0 > t' >t$.
Since $P$ is finitely generated, using Nakayama's lemma, we can conclude that $\phi^{t'}_t$ is an isomorphism.

\begin{claim}\label{CL:E-FINITE}
$I/E$ is finite
\end{claim}
Claim~\ref{CL:OPEN} implies that each equivalence class of $E$ has
non-empty interior. Since $E$ is definable, $I/E$ must be finite.

It follows from Claim~\ref{CL:E-FINITE} that  there exists a left
neighbourhood $J$ of $0$ such that, for every $t,t' \in J$, $t E t'$.
Therefore, for every $t \leq t'$ in that neighbourhood $\phi^{t'}_t$
is an isomorphism, and we conclude by Remark~\ref{lim}.
\end{proof}

\br
In the situation of Theorem~\ref{THM:DEFINABLE}, let
\begin{align*}
f &: Y \to M^{\geq 0}\\
f(x) &:= \inf \set{t: x \in Y_t}.
\end{align*}
We have remarked that $f$ might not be continuous.
Assume that there exists $t_0 >0$ such that $f$ is continuous on $Y_{t_0}$.
Then, the proof of Theorem~\ref{THM:DEFINABLE} can be simplified.
In fact, by the trivialization theorem, we can assume that there exists a continuous definable map $\lambda: Y_{t_0} \setminus A \to F$ (where $F := f^{-1}(t_0)$), such that the map
\begin{align*}
\mu:= (f,\lambda) &: Y_{t_0}\setminus A \to (0,t_0) \times F
\end{align*}
is a homeomorphism.
Let $0 < t \leq t' \leq t_0$.
Let $\theta: [0,1] \times [0,t'] \to [0,t']$ be a definable 
strong deformation retraction between $[0,t']$ and $[0,t]$.
Define
\[\begin{array}{r@{\,}ll}
\Lambda: [0,1] \times Y_{t'} &\to Y_{t'}\\
\Lambda\Pa{s, \mu^{-1}(u,x)} &= \mu^{-1}\Pa{\theta(s,u),x} &\text{if }  0 <u <t',\\
\Lambda(s,a) &= a &\text{if } a \in A.
\end{array}\]
Note that $\Lambda$ gives a definable strong deformation retraction between
$Y_t'$ and~$Y_t$.
Therefore the inclusion $Y_t \subseteq Y_{t'}$ induces an isomorphism in cohomology.
\er
The fact that in the situation of the above remark $Y_t$ is a deformation retract of $Y_{t'}$ follows from a stronger results in \cite[Thm.~A.5]{PeterzilS07}, where however it is assumed that $0$ is a regular value for $f$.

\bq Under the hypothesis of Theorem~\ref{THM:DEFINABLE}, is there some
$t_0 >0$ such that $A$ is a deformation retract of $Y_t$ for every $0
< t < t_0$?  \eq

Note that in the proof of Theorem~\ref{THM:DEFINABLE}
we have used heavily the trivialization theorem, which holds only for
expansions of fields.  

\bq Can we weaken the hypothesis of
Theorem~\ref{THM:DEFINABLE} dropping the condition that $M$ expands a
field \rom(but saying instead that it expands a group\rom)?  \eq

\appendix
\appendixpage
\addappheadtotoc

\section{Derived functors}
Let $\Acat$ and $\Bcat$ be Abelian categories such that $\Acat$ has
enough injective objects (namely, every object of $\Acat$ can be
embedded in an injective one), and let $\Gamma : \Acat \to \Bcat$ be
an additive covariant left-exact functor.

In our applications, $\Acat$ will be the category $\Sh_X$ of sheaves
of Abelian groups on a fixed topological space $X$, $\Bcat$ the
category of Abelian groups $\Ab$, and $\Gamma \colon \Sh_X \to \Ab$
the global section functor.

A {\bf chain complex} in $\Acat$ is a sequence of maps $C^n \stackrel {d^n}
\to C^{n+1}$ in $\Acat$ with $d^{n+1} \circ d^{n} = 0$ for every $n
\in \Z$. The chain complex $C^* = (C^n,d^n)_{n\in \Z}$ is {\bf exact} if
$\Ker(d^n)=\im(d^{n-1})$ for every $n$. 
A {\bf resolution} of $A\in \Acat$ is an exact chain complex of the form
\[
0 \to A \stackrel {j} \to I^0 \stackrel {d^0} \to I^1
\stackrel {d^1} \to I^2 \dotso
\]

Given such a resolution, after applying $\Gamma$
we obtain a chain complex in~$\Bcat$ 
\begin{equation} 
0 \to \Gamma(A) \stackrel {\Gamma(j)} \to \Gamma(I^0)  \stackrel {\Gamma(d^0)} \to \Gamma(I^1) \stackrel {\Gamma(d^1)} \to \Gamma(I^2) 
\ldots
\end{equation} 
which however need not be exact since $\Gamma$ is only left-exact. The lack of exactness is measured by 
the cohomology groups 
\[\Hom^n(\Gamma(I^*)) =\Ker(\Gamma(d^n))/\im(\Gamma(d^{n-1}))\] 
where we stipulate that 
$d^{-1}:= 0$. Note that $\Hom^0(\Gamma(I^*)) = \Gamma(A)$. 

Now let us further assume that $0\to A \to I^*$ is an {\bf injective
resolution}, namely a resolution such that each $I^n$ is an injective
object of $\Acat$.  Then it turns out that $\Hom^*(\Gamma I^*)$ depends only on $A$ and not on the particular choice of the injective resolution. This can be seen as follows. 

If $0 \to B \to C^*$ is another exact sequences in
$\Acat$ and $f \colon B \to A$ is a morphism in $\Acat$, then (by
injectivity of $I^n$) there exists a morphism of chain complexes $
f^* \colon C^* \to I^* $ commuting with $f$ as in Diagram 1 below.

\begin{center}
\begin{minipage}{\textwidth}
\DIAGV{50} {0} \n{\Ear{}} \n{B} \n{\Ear{}} \n{C^*} \nn 
           {}\n{} \n{\Sar{f}}   \n{}     \n{\Sar{f^*}} \nn 
           {0}\n{\Ear{}} \n{A} \n{\Ear{}} \n{I^*}        
\diag\vspace{-2.5ex}
\centerline{Diagram 1}
\end{minipage}
\end{center} 

Moreover $f^*$ is determined by
$f$ up to chain homotopy (see \cite[Proposition
A3.13]{EISENBUD:1995}). Since $\Gamma$ preserves chain homotopies, we
obtain a morphism (unique up to chain homotopy) \beq \Gamma(f^*) \colon
\Gamma(C^*) \to \Gamma(I^*)\eeq which induces a unique map \beq
\Hom^n(\Gamma(f^*))\colon \Hom^n(\Gamma(C^*))\to \Hom^n(\Gamma(I^*)). \eeq It
follows that if $0\to A \to I^*$ and $0 \to A \to C^*$ are two
injective resolutions of $A \in \Acat$, then the identity morphism
$\id_A$ on $A$ induces a canonical isomorphism \beq
\Hom^n(\Gamma(\id_A))\colon \Hom^n (\Gamma(I^*)) \cong \Hom^n (\Gamma(C^*)). \eeq

The right derived functor $\der n \Gamma \colon \Acat \to \Bcat$ can thus be defined, up to a canonical isomorphism, by:
\beq
\der n \Gamma (A) := \Hom^n(\Gamma I^*)
\eeq
If $f$ is as in the Diagram 1 and all the $I^n,C^n$ are acyclic, we define
\beq
\der n \Gamma (f) := \Hom^n(\Gamma f^*)
\eeq

An object $C \in \Acat$ is {\bf acyclic} (for
$\Gamma$), if $\der n\Gamma (C) = 0$ for all $n>0$. Any injective
object in $\Acat$ is acyclic for any additive functor. 
If $A \in \Acat$ and $0\to A \to C^*$
is an acyclic resolution of $A$, then the natural map
$\Hom^n(\Gamma(C^*)) \to \der n \Gamma (A)$
is an isomorphism (the proof in \cite[\S{}II.4]{Bredon97}
works in this context). 
So we can compute $\der n \Gamma (A)$ using any acyclic resolution of~$A$, 
not necessarily injective. 
We have the following sufficient condition for acyclicity:

\bl\label{LEM:ACYCLIC}
Let $\Ccat$ be a class of objects in~$\Acat$.
Assume that
\begin{enumerate}
\item All injective objects of $\Acat$ are in~$\Ccat$\rom;
\item If $0 \to C' \to A \to A'' \to 0$ is a short exact sequence \rom(in $\Acat$\rom), with $C' \in \Ccat$, then
\[
0 \to \Gamma(C') \to \Gamma(A) \to \Gamma(A'') \to 0
\]
is a short exact sequence \rom(in $\Bcat$\rom)\rom;
\footnote{Note that, since $\Gamma$ is left-exact, $0 \to \Gamma(C') \to \Gamma(A) \to \Gamma(A'')$ is always exact.}
\item If $0 \to C' \to C \to A'' \to 0$ is a short exact sequence, with $C', C \in \Ccat$, then $A'' \in \Ccat$.
\end{enumerate}
Then,
\begin{description}
\item[(a)] If $0 \to C^{-1} \to C^0 \to C^1 \dotso$ is an exact sequence in $\Ccat$, then $0 \to \Gamma(C^{-1}) \to \Gamma(C^0) \to \Gamma(C^1) \dotso$ is also exact\rom;
\item[(b)] $\der n \Gamma(C) = 0$ for every $C \in \Ccat$ and $n > 0$ \rom(namely, all objects in $\Ccat$ are acyclic\rom).
\end{description}
\el
\begin{proof} Point (b) follows from (a).
The proof of (a) is in \cite[Thm.~II.3.1.3]{Godement73}
for the class $\Ccat$ of flabby sheaves: it generalizes easily to this context. 
\end{proof}

\section{Flabby and soft sheaves}
Let $X$ be an arbitrary topological space (not necessarily Hausdorff).
We shall give some results about sheaf cohomology on $X$; for unproved facts in this and the following sections, the reader can look either in \cite{Bredon97} or in \cite{Godement73}.
\bd\label{DEF:COM}
We shall denote by $\Xcat$ the category of sheaves of Abelian groups on~$X$.
$\Xcat$ is an Abelian category, with enough injective objects.
Let
\[\begin{array}{r@{\ }r@{\,}c@{\,}l}
\Gamma: & \Xcat  & \to     & \Ab\\
        & \ca F  & \mapsto & \ca F(X)
\end{array}\]
be the global section functor.
It is well-known that $\Gamma$ is left-exact.
The $n$-th cohomology functor on $X$ is the $n$-th right derived functor of~$\Gamma$
\[
\Hom^n(X; \ca F) := \der n \Gamma (\ca F).
\]
Given $A \subseteq X$ we shall write $\Hom^*(A; \ca F)$ for $\Hom^*(A; \ca F \rest A)$.
\ed
A sheaf on $X$ is called {\bf flabby} if for every open subset $A$ of
$X$ any section of ${\ca F}$ on $A$ can be extended to the whole of
$X$. Any injective sheaf is flabby, and any flabby sheaf is acyclic.
A sheaf $\ca F$ on $X$ is {\bf soft} if for
every closed subset $A$ of $X$ any section of $\ca F$ on $A$ can be
extended to the whole of $X$. On Hausdorff paracompact spaces, soft
sheaves are flabby, and therefore acyclic.

Many results on sheaf cohomology
are stated in \cite{Godement73} or \cite{Bredon97} for paracompact
Hausdorff spaces, but an analysis of the proofs shows that they hold
more generally for arbitrary \PCN spaces (cf.\ \S\ref{preliminaries}).

\bl\label{LEM:PCN-EXT} Let $X$ be a topological space and let $A
\subseteq X$ have a basis of \PCN neighbourhoods.  Let ${\ca F} \in
\Xcat$, and $s$ be a section of ${\ca F}$ on~$A$.  Then, $s$ can be
extended to a neighbourhood of $A$.  Therefore,
\[
{\ca F}(A) = \dirlim_{A \subseteq U} {\ca F}(U),
\]
using the neighbourhoods $U$ of~$A$, ordered by reversed inclusion.
\el

\bp
The proof of \cite[Thm. 3.3.1, p. 150]{Godement73} works under our hypothesis.
\ep

\bl\label{LEM:SOFT} Let ${\ca F}$ be a flabby sheaf on $X$.  If $X$ is
\PCN, then ${\ca F}$ is soft.  More generally, if $X$ is an arbitrary
topological space, ${\ca F}$ is a flabby sheaf on $X$, and $A$ is a subspace of
$X$ with a basis of \PCN neighbourhoods, then ${\ca F}$ induces a soft
sheaf on $A$. \el
\begin{proof}
Let $C \subseteq A$ be closed in~$A$, and $s$ be a section of ${\ca F}$ on~$C$.
By Lemma~\ref{LEM:PCN-EXT}, $s$~can be extended to a neighbourhood $W$ of $C$ in~$X$.
Since ${\ca F}$ is flabby, $s$~can be further extended to~$X$, and \emph{a fortiori} to~$A$.
\end{proof}

\bl\label{REM:CEMBEDDED} If $X$ is a topological space and $A
\subseteq X$ has a basis of \PCN neighbourhoods in~$X$, then $A$ is
\PCN.  More precisely, if $A \subseteq X$ has a basis of normal
neighbourhoods, then $A$ is normal.  If $A$ has a basis of paracompact
neighbourhoods, then $A$ is paracompact.  \el

\bp Normality: First note that a space is normal if and only if every
open covering consisting of two open sets $U,V$ has a shrinking. Given
an open covering $\{U,V\}$ of $A$, let $U',V'$ be open sets of $X$
with $U'\cap A = U$ and $V'\cap A = V$. Let $W$ be a normal open
neighbourhood of $A$ contained in $U'\cup V'$. We can assume $W =
U'\cup V'$. By the normality of $W$ we can find a shrinking $\{U'_0,
V'_0\}$ of $\{U', V'\}$. The intersections of $U'_0$ and $V'_0$ with
$A$ is a shrinking of $\{U,V\}$.  The proof of paracompactness is
similar. 
\ep 

\bt \label{LEM:PCN-ACYCLIC}
Assume that $X$ is a \PCN space. Then any soft sheaf on $X$ is acyclic. 
\et
\begin{proof}
(See \cite[\S{}II-3.5]{Godement73}.)  Lemma~\ref{LEM:ACYCLIC} gives a
sufficient condition for acyclicity, so it suffices to verify its
three hypothesis. 

Point 1. We need to show that any injective sheaf on $X$ is soft. This
follows from the fact that any injective sheaf is flabby, and any
flabby sheaf on a \PCN space is soft (Lemma~\ref{LEM:SOFT}).

Point 2. This is proved in~\cite[Thm.~II-3.5.2, p. 153]{Godement73}
for Hausdorff paracompact spaces (or supports), but the proof works
also for \PCN spaces.

Point 3. Assume that $0 \to {\ca F}' \to {\ca F} \to {\ca F}'' \to 0$
is an exact sequence in~$\Xcat$, with ${\ca F}'$ and ${\ca F}$ soft.  We
have to prove that ${\ca F}''$ is also soft.  Let $C \subseteq X$ be
closed.  By the second point, since ${\ca F}' \rest C$ is soft,
\[
{\ca F}(C) \to {\ca F}''(C)
\]
is surjective.
Hence, a section $s''$ of ${\ca F}''$ on $C$ is represented by a section $s$ of ${\ca F}$ on~$C$.
Since ${\ca F}$ is soft, $s$ can be extended to all~$X$, and hence also $s''$ can be extended to~$X$.
\end{proof}

\section{Taut subspaces}

\bd Let $A \subseteq X$, and ${\ca F} \in \Xcat$.  For every neighbourhood
$U$ of $A$ in $X$, the inclusion map $\lambda^A_U: A \to U$
induces a map in cohomology $\Hom^*(\lambda^A_U): \Hom^*(X;{\ca F}) \to
\Hom^*(A; {\ca F})$. 
If $U \subseteq V$, the maps
$\Hom^*(\lambda^U_V)$, $\Hom^*(\lambda^A_U)$, and
$\Hom^*(\lambda^A_V)$ commute.  Hence, we have a canonical map
\[
\dirlim_{A \subseteq U} \Hom^*(U; {\ca F}) \to \Hom^*(A; {\ca F} ).
\]
We say that $A$ is {\bf taut in $X$ for the sheaf ${\ca F}$} if the above
map is an isomorphism.  We say that $A$ is {\bf taut in $X$} iff $A$ is taut
for every sheaf ${\ca F}$.
\ed

We give an example of a normal \PCN space $X$ with a
quasi-compact subspace $A$ which is not taut in $X$.

\bexa Let $X = \{a,b,c,d\}$ partially ordered so that $a = \min X, d =
\max X$ and $b,c$ are incomparable. Put the following topology on $X$:
a set is closed iff it is downward closed in this order. Then $X$ is a
normal spectral space. The subset $A = \{b,c\}$ of $X$ is
quasi-compact (since it is finite) but not taut. In fact $V =
\{b,c,d\}$ is the smallest open set containing $A$, but $V$ is
connected and $A$ is disconnected. Hence 
\[
\Hom^0(A, \Z) = \Z^2 \neq \Z = \Hom^0(V, \Z) = \dirlim_{A \subseteq U} \Hom^0(U, \Z).
\]
\eexa

The above example shows that if we want to prove that a certain subspace is
taut, we need some normality assumptions on its neighbourhoods. 

\bl
\label{LEM:EMBEDDED-FAMILY}
Let $A \subseteq X$, and $\Pa{Y_i}_{i \in I}$ be a family of subsets
of $X$ indexed by a filtered set $I$, such that $\Pa{Y_i}_{i \in I}$
is decreasing, and $A \subseteq \bigcap_{i \in I} Y_i$.  Assume that
$A$ and each $Y_i$ are taut in $X$, and that, for every neighbourhood
$V$ of $A$, there exists $i \in I$, such that $Y_i \subseteq V$.
Then for every ${\ca F} \in \Xcat$,
\[
\Hom^*(A; {\ca F}) 
= \dirlim_{i \in I} \Hom^*(Y_i; {\ca F}). 
\]
\el
\begin{proof}
Since $A$ and each $Y_i$ are taut in $X$, 
\[
\Hom^*(A; {\ca F}) 
= \dirlim_{A \subseteq U} \Hom^*(U; {\ca F}) =
\dirlim_{i \in I} \dirlim_{Y_i \subseteq U} \Hom^*(U; {\ca F}) =
\dirlim_{i \in I} \Hom^*(Y_i; {\ca F} 
). \qedhere
\]
\end{proof}

\bt\label{LEM:TAUT-PCN}
Let $X$ be an arbitrary topological space and let $A \subseteq X$.
Assume that $A$ has a basis of \PCN neighbourhoods in~$X$.
Then $A$ is taut in~$X$. 
\et
\bp 
(See \cite[Thm.~II.4.11.1]{Godement73}.)
Let ${\ca F} \in \Xcat$, and $0 \to {\ca F} \to {\ca C}^*$ be an injective resolution (in~$\Xcat$). 
By Lemma~\ref{LEM:PCN-EXT}, ${\ca C}^*(A) = \dirlim_{A \in U} {\ca C}^*(U)$.
Since any injective sheaf is flabby (see \cite{Bredon97}),
by Lemma~\ref{LEM:SOFT},
\[
0 \to {\ca F} \rest A \to {\ca C}^* \rest A
\]
is a soft resolution in $\Sh_A$. Moreover $A$ is \PCN (Lemma
\ref{REM:CEMBEDDED}), so any soft sheaf on $A$ is acyclic (Lemma
\ref{LEM:PCN-ACYCLIC}) and therefore the cohomology of ${\ca F} \rest A$
can be computed by soft resolutions. So we obtain:
\[
\Hom^p(A; {\ca F} \rest A) = \Hom^p\Pa{\Gamma_A({\ca C}^*\rest A)}
= \Hom^p\Pa{\dirlim_{A \subseteq U} {\ca C}^*(U)}
= \dirlim_{A \subseteq U} \Hom^p(U; {\ca F}). \qedhere
\]
where the last equality follows from the fact that the homology
functor from complexes of modules over a fixed ring~$R$, to
$R$-modules, preserves inductive limits over filtered sets (see
\cite[\S I.2.1]{Godement73}). 
\ep 

\bc\label{taut}
{\rm(\cite{Delfs85,Jones06,EdmundoJP05})}  Let $X$ be a definable set in
an o-minimal expansion of a group. Let $A$ be a quasi-compact subset
of $\widetilde X$ \rom(this assumption holds in particular if $A$ is
definable, or type-definable\rom).  Then $A$ is taut in~$\widetilde
X$. \ec

\bc\label{LEM:TAUT-DEFINABLE} Let $\Pa{Y_t}_{t > 0}$ be a definable
family of definably compact subsets of some definable set~$Y$, such
that $Y_{t'} \subseteq Y_t$ for every $0 < t' < t$.  Let $A :=
\bigcap_{t > 0} Y_t$.  Then the inclusion induces an isomorphism
\[
\Hom^*\Pa{ A; {\ca F}} \cong \dirlim_{t \to 0} \Hom^*\Pa{{Y_t}; {\ca F}},
\]
for every sheaf  ${\ca F}$ on~$\widetilde Y$.
\ec
\begin{proof}
Since the $Y_t$ are definably compact, for every definable
neighbourhood $U$ of $A$ there exists $t >0$ such that $Y_t \subseteq U$.
Moreover $A$ and $Y_t$ are taut by Theorem \ref{LEM:TAUT-PCN} or
Corollary \ref{taut}. Hence we can apply
Lemma~\ref{LEM:EMBEDDED-FAMILY}
\end{proof}

\section{\Cech cohomology}\label{AP:D}
\begin{definition}
For every open covering $\Ucov$ of~$X$ and sheaf $\Fsh$ on $X$, we can define the \Cech cohomology groups $\HomC^*(\Ucov;\Fsh)$ as in~\cite{Bredon97}.
If $\Vcov$ is a refinement of $\Ucov$, there exists a canonical map $\HomC^*(\Ucov; \Fsh) \to \HomC^*(\Vcov; \Fsh)$.
We remind that the \Cech cohomology groups of $X$ are defined as
$
\HomC^*(X; \Fsh) := \dirlim_{\Ucov} \HomC(\Ucov; \Fsh).
$
\end{definition}
\begin{lemma}\label{LEM:LOC-0}
Let $X$ be a \PCN space and $\ca F$ be a presheaf on~$X$, such that 
$\Fshc = 0$, where $\Fshc$ is the sheaf generated by~$\ca F$.
Then, $\HomC^*(X;\ca F) = 0$.
\end{lemma}
\begin{proof}
Bredon~\cite[Thm.~III.4.4]{Bredon97}, Godement~\cite[Thm.~II.5.10.2]{Godement73} and Spanier~\cite[Thm.~6.7.16]{Spanier81} prove the above lemma for the case when $X$ is Hausdorff and paracompact.
Their proofs work also for $X$ \PCN, without need of modifications.
\end{proof}

\begin{lemma}\label{LEM:GENERATED}
Let $X$ be a \PCN space, and $\ca F$ be a presheaf on~$X$.
Denote by $\Fshc$ the sheaf generated by~$\ca F$.
Then, the canonical map $\theta: \ca F \to \Fshc$ induces an isomorphism
\[
\HomC^*(X; \ca F) 
\Ar{\cong}\HomC^*(X;\Fshc).
\]
\end{lemma}
\begin{lemma}\label{LEM:CECH}
Let $X$ be a \PCN space, and $\ca F$ be a sheaf on it.
Then, there is a natural isomorphism
\[
\HomC^*(X;\ca F) \cong \Hom^*(X;\ca F).
\]
\end{lemma}
\begin{proof}[Proof of Lemmata~\ref{LEM:GENERATED} and~\ref{LEM:CECH}]
Bredon~\cite[Cor.~III.4.5, Cor.~III.4.12]{Bredon97} and Godement~\cite[Cor.~II.5.10]{Godement73} prove the above lemmata for the case when $X$ is Hausdorff and paracompact.
The proof for the \PCN case can be done as in~\cite{Bredon97}.
The main  ingredient is Lemma~\ref{LEM:LOC-0};
the rest are algebraic manipulations that do not use any property of the space.
\end{proof}
In the case of $X$ normal spectral space, Lemma~\ref{LEM:CECH} was already proven in~\cite{CarralC83}.

\medskip

\begin{footnotesize}
\noindent {\sc Universit\`a di Pisa, Dipartimento di Matematica, Largo
Bruno Pontecorvo 5, 56127 Pisa, Italy.} E-mail: {\tt berardu@dm.unipi.it}.
\medskip

\noindent {\sc Universit\`a di Pisa, Dipartimento di Matematica, Largo
Bruno Pontecorvo 5, 56127 Pisa, Italy}. E-mail:
{\tt fornasiero@dm.mail.unipi.it}.
\end{footnotesize} 


\begin{thebibliography}{CosteR82} 

\bibitem[Bred97]{Bredon97} G. E. Bredon. Sheaf theory. Graduate Texts in Mathematics (Vol. 170),
Second edition. Springer-Verlag, New York, 1997. xii+502 pp.

\bibitem[CarrC83]{CarralC83} M. Carral, M. Coste. Normal spectral spaces and
their dimension, J. Pure Appl. Algebra 30 (1983) 227-235

\bibitem[CosteR82]{CosteR82} M. Coste, M.-F. Roy, La topologie du Spectre Reel,
in: D. W. Dubois and T. Recio, eds., Ordered Fields and Real Algebraic
Geometry, Contemporary Mathematics 8, 1982, pp. 27-59

\bibitem[Delfs85]{Delfs85} H. Delfs. The homotopy axiom in semialgebraic
cohomology. Journal f\"ur die reine und angewandte Mathematik 355
(1985) 108-128

\bibitem[Dries98]{Dries98} L. van den Dries. Tame Topology and o--minimal
structures. London Math. Soc. Lecture Notes Series, vol. 248. Cambridge
Univ. Press 1998. x+180 pp.

\bibitem[God73]{Godement73} Roger Godement.
Topologie alg\'ebrique et th\'eorie des faisceaux.
Troisi\`me \'edition. Publications de l'Institut de
Math\`ematique de l'Universit\'e de Strasbourg, XIII. Actualit\'s
Scientifiques et Industrielles, No. 1252. Hermann, Paris,
1973. viii+283 pp.

\bibitem[EdmJP05]{EdmundoJP05} M. Edmundo, G. O. Jones, N. J. Peatfield, Sheaf
cohomology in o-minimal structures, Journal of Mathematical Logic, Vol. 6, No. 2 (2006) 163 179

\bibitem[EilS52]{EilenbergS52} S. Eilenberg, N. Steenrod. Foundations
of algebraic topology. Princeton University Press, Princeton, New
Jersey, 1952. xv+328 pp.

\bibitem[Eisen95]{EISENBUD:1995}
D. Eisenbud. Commutative algebra, With a view toward algebraic geometry. Graduate Texts in
  Mathematics (vol. 150). Springer-Verlag, New York, 1995.

\bibitem[Eng89]{Engelking89} R. Engelking. General topology. Rev. and compl. ed. (English)
Sigma Series in Pure Mathematics, 6. Berlin: Heldermann Verlag. viii, 529 p. 

\bibitem[Jones06]{Jones06} G. O. Jones. Local to global methods in o-minimal
expansions of fields, PhD thesis, Oxford 2006, pp. 83

\bibitem[PetSS00]{PeterzilSS00} Y. Peterzil, P. Speissegger, S. Starchenko.
Adding multiplication to an o-minimal expansion of the additive group of real numbers.  
Logic Colloquium '98 (Prague),  Lecture Notes in Logic 13, 357-362. Assoc. Symbol. Logic, Urbana, IL, 2000. 

\bibitem[PetS07]{PeterzilS07}
Y. Peterzil, S. Starchenko. 
Computing o-minimal topological invariants using differential topology. 
Trans. Amer. Math. Soc. 359 (2007), no. 3, 1375--1401

\bibitem[Pillay87]{Pillay87}
A. Pillay. First order topological structures and theories. The Journal of
Symbolic Logic, vol. 52, n. 3 (1987) 763 - 778

\bibitem[Pillay88]{Pillay88}
A. Pillay, Sheaves of continuous definable functions. The Journal of
Symbolic Logic, vol. 53, n. 4 (1988) 1165 - 1169

\bibitem[Span81]{Spanier81}
E.~H. Spanier. Algebraic topology.
Springer-Verlag, New York, 1981.
Corrected reprint of the 1966 original.

\bibitem[Wallm38]{Wallman38}
H. Wallman. Lattices and Topological Spaces. Annals of Mathematics (2) 39 (1938) 112-126
\end{thebibliography}
\end{document}